\documentclass[trsc,nonblindrev]{informs3} 

\OneAndAHalfSpacedXI 

\usepackage{xcolor}
\usepackage{colortbl,booktabs}
\usepackage{booktabs}
\usepackage{multirow} 
\usepackage[figuresright]{rotating} 
\usepackage{adjustbox}
\usepackage{pdflscape}

\usepackage{verbatim}

\usepackage{algorithm}  
\usepackage{algorithmic}   

\usepackage{natbib}
\bibpunct[, ]{(}{)}{,}{a}{}{,}%
\def\bibfont{\small}%

\usepackage{amsmath} 

\TheoremsNumberedThrough     

\EquationsNumberedThrough    

\usepackage{cleveref}
\crefname{section}{§}{§§}
\Crefname{section}{§}{§§}

\begin{document}
	
	
	\RUNAUTHOR{Wang et al.}
	
	\RUNTITLE{Tactical fleet and crew scheduling}
	
	\TITLE{Tactical fleet assignment and crew pairing problem with crew flight time allocation}
	
	\ARTICLEAUTHORS{%
		\AUTHOR{Danni Wang, Siqi Guo, Wenshu Wang, Zhe Liang}
		\AFF{School of Economics and Management, Tongji University, \EMAIL{wangdanni@tongji.edu.cn},
		\EMAIL{1930344@tongji.edu.cn},
		\EMAIL{2010083@tongji.edu.cn}, \EMAIL{liangzhe@tongji.edu.cn}}
	} 
	
	\ABSTRACT{%
		Aircraft and crew are two major resources that ensure the smooth operations of airlines. However, with the anticipated growth in the aviation industry, the crew resource is predicted to be insufficient worldwide and has been one of the bottlenecks in fast-developing airlines. In addition, the mismatch between aircraft and crew has influenced the airlines' operation, and limited the full usage of all resources. To resolve this problem, we propose a tactical fleet assignment and crew pairing problem with crew flight time allocation (TFACPP). The basic integrated model is reformulated by the Benders decomposition, where the Benders master problem (BMP) poses the most significant computational barrier. To efficiently solve the BMP, we propose a column generation algorithm. The TFACPP can provide superior solutions compared to the equal allocation of crew resources. In addition, we provide a quantitative method for evaluating the scarcities of crew and aircraft resources and the matching degree between crew and aircraft based on the shadow prices of the proposed model. These information can provide rich managerial insights regarding the acquisition, replacement, and transition of crew and aircraft. 
		
	}%
	
	
	\KEYWORDS{tactical fleet and crew pairing problem, crew flight time allocation}
	
	\maketitle
	
	%
	

\section{Introduction}\label{intro}	
The aviation industry has become a significant player in the global
economy, and people increasingly prefer to air travel. The International Civil Aviation Organization (ICAO) anticipated that the global passenger traffic will grow at 4.1\% annually from 2015 to 2045 \citep{ICAO2018}. The revenue passenger kilometer (RPK) value in 2018 was 8,255 million, marking the ninth consecutive year of above-trend RPK growth. To accommodate the market, airlines are aggressively expanding their fleet. The International Air Transport Association (IATA) reported that the fleet increased by more than 1,000 aircraft (to 29,754) in 2018 compared to 2017 \citep{IATA2018, IATA2019}. Although the aviation industry was negatively affected by the COVID-19 pandemic at the beginning of 2020, the long-term development trend remains positive. According to the Civil Aviation Administration of China (CAAC), the passenger throughput of Chinese airports reached 857.159 million in 2020, which had recovered to the level before the COVID-19 pandemic (2019) of 63.2\% \citep{CAAC2021}.

The expansion of airlines poses a tremendous challenge for aviation scheduling and recovery \citep{Belobaba2009,Barnhart2012,Eltoukhy2017,Zhou2020,Su2021}. The airline planning is generally divided into four sequential subproblems, namely, schedule design problem (SDP), fleet assignment problem (FAP), aircraft routing problem (ARP), and crew scheduling problem (CSP). The SDP involves developing a timetable that includes the origin and destination as well as the departure and arrival time of each flight. The FAP involves determining the fleet type of aircraft to operate each flight while maximizing the total profit. The ARP involves assigning an aircraft for each flight while satisfying the maintenance requirements. Given the solution of the ARP, the CSP is solved to cover the scheduled flights by the available crews without violating relevant rules \citep{Stojkovic1998}. Although significant progress has been achieved in the above sequential approaches, many researchers and practitioners are still devoting considerable efforts to the new emerging problems in the aviation industry. This is not only due to theoretical and technical obstacles, but also the continuously changing and increasingly restrictive business environments and configurations. Next, we present three challenges in current practice and highlight the main contributions of our paper.

First, the aggravated mismatch between the crew and aircraft resources has affected airline operations seriously. Because pilot development is expensive and time-consuming, the pilot supply has continued to decrease worldwide \citep{AINonline2017,CNN2018,Statista2022}. Boeing forecast that 602,000 new pilots will be needed over the next 20 years, with 134,000 and 120,0000 new pilot requirements in North America and European, respectively \citep{Boeing2022a}. The situation is worse in rapidly developing markets such as China, where the pilot supply is always insufficient with the quickly expanded fleet. It is predicted that Chinese airlines would require 8,485 new airplanes over the next 20 years, or an average of more than 400 aircraft each year by 2041 \citep{Boeing2022b}. This will definitely exacerbate the mismatch between the crew and aircraft resources. In practice, it is typical for the flight schedule, which is designed based on aircraft availability, to be canceled due to a lack of crew resources.

Second, studies considering the macroscopic crew resource availability are rare. Specifically, in the existing literature, the CSP primarily focuses on operational rules, such as the 8-in-24 rule \citep{Barnhart2003,Kasirzadeh2015}. However, as shown in Table \ref{tab_limits}, there are various rules limiting the crew flight time over different time periods. For instance, FAA restricts that a pilot can fly up to 8-9 hours in a working day, whereas only 4 hours in average if we consider the 30-day time horizon. This number further decreases to 3.3 and 2.7 hours per day in the 90-day and 12-month time horizon, respectively. It is easy to notice that the macroscopic rules are more stringent than the microscopic rules. The situation is similar in Europe and China. As a result, the yearly crew flight time is commonly one of the bottlenecks of the CSP in practice. However, this difficulty is often overlooked because the researchers focus more on short-term schedules.

\begin{table}[]
	\caption{Flight time limits over a rolling time period}
	\centering
	\small
	{\def\arraystretch{1.25}
		\begin{tabular}{lrlrl}
			\toprule
			\multicolumn{1}{c}{\begin{tabular}[c]{@{}c@{}}Admini-\\ stration\end{tabular}} & \multicolumn{1}{c}{\begin{tabular}[c]{@{}c@{}}Flight time \\ limit (hours)\end{tabular}} & \multicolumn{1}{c}{\begin{tabular}[c]{@{}c@{}}Time \\ horizon\end{tabular}} & \multicolumn{1}{c}{\begin{tabular}[c]{@{}c@{}}Average flight time\\ per day (hour)\end{tabular}} & \multicolumn{1}{c}{\begin{tabular}[c]{@{}c@{}}\\ Source\end{tabular}}  \\
			\midrule
			FAA                                & 8-9                                                                                       & 1 day                                                & 8-9                                                                                                & CFR.14.117.11 \citep{FAA2013}              \\
			FAA                                & 120                                                                                       & 30 days                                               & 4                                                                                                  & CFR.14.121.483 \citep{FAA1996}             \\
			FAA                                & 300                                                                                       & 90 days                                              & 3.3                                                                                                & CFR.14.121.483 \citep{FAA1996}             \\
			FAA                                & 1000                                                                                      & 12 months                                                 & 2.7                                                                                                & CFR.14.121.483 \citep{FAA1996}             \\
			EASA                               & 100                                                                                       & 28 days                                                & 3.5                                                                                                & ORO.FTL.210 \citep{EASA2019}               \\
			EASA                               & 1000                                                                                      & 12 months                                                 & 2.7                                                                                                & ORO.FTL.210 \citep{EASA2019}               \\
			EASA                               & 900                                                                                       & 1 calendar year                                                & 2.5                                                                                                & ORO.FTL.210 \citep{EASA2019}               \\
			CAAC                               & 8-9                                                                                       & 1 day                                                & 8-9                                                                                                & CCAR.121.483 \citep{CAAC2017}              \\
			CAAC                               & 100                                                                                       & 1 month                                                 & 3.3                                                                                                & CCAR.121.487 \citep{CAAC2017}               \\
			CAAC                               & 900                                                                                       & 1 calendar year                                                 & 2.5                                                                                                & CCAR.121.487 \citep{CAAC2017}     \\
			\bottomrule         
		\end{tabular}
	}
	\label{tab_limits}%
\end{table}

Finally, airlines need to allocate the yearly crew flight time into different months properly according to demand. Currently, the industrial allocating approach primarily depends on experts' experience. In practice, it is quite common for airlines to equally distribute the yearly crew flight time into each month in the planning stage, and then make adjustments for the peak seasons. However, this may not capture the mild month-to-month fluctuations in demand. In addition, with the spread of COVID-19, demand fluctuations have grown more frequent and severe. Thus, there is an urgent need to develop a systemic method to allocate the yearly crew flight time into different months to match the changing seasonal demand.

To overcome the above-mentioned challenges, we address a tactical fleet assignment and crew pairing problem with crew flight time allocation (TFACPP), and formulate it as a large-scale mixed integer programming (MIP) model. The TFACPP simultaneously optimizes the fleet assignment and crew pairing problems in consideration of macroscopical crew flight time constraints. Moreover, we extend the integrated model to include crew transition from one fleet family to another and uncertainty of crew availability. The crew transition is frequently used in the industry but is seldom investigated by academics. The uncertainty is caused by unexpected personnel absences, such as illness and vacation.

Another contribution of our paper lies in the solution method. We first reformulate the basic TFACPP model by the Benders decomposition. Furthermore, we demonstrate that the principal solution challenge pertains to solving the Benders master problem (BMP). Thus, we design a column generation algorithm to solve the BMP.

Finally, we demonstrate that the dual variables for the aircraft (crew) availability constraints can provide rich managerial information for long-term decisions, such as aircraft replacement and acquisition and crew training planning. Currently, such decisions rely heavily on the expertise of airline managers due to the lack of a quantitative method. Given the significance of these decisions, airlines are urgently seeking comprehensive methods to evaluate the scarcities of different types of resources and supplement them as appropriate.

The remaining paper is organized as follows. Section \ref{review} presents a review of the existing work related to the TFACPP. Section \ref{model} describes the formulation of the MIP model for the TFACPP. Section \ref{method} introduces the Benders reformulation and solution method. Section \ref{extend} presents two types of extensions pertaining to industrial constraints. Several computational experiments are discussed in Section \ref{test}, and practical managerial insights are presented in Section \ref{insight}. Finally, we provide our conclusions and future work in Section \ref{conclu}.

\section{Literature review}\label{review}
This section reviews the fleet assignment problem and crew pairing problem, followed by the integrated fleet assignment and crew pairing problem.

\subsection{Fleet assignment problem} \label{review_fap}
A typical FAP covers the scheduled flights with different fleet types to maximize profit. There are two classic fleet assignment models (FAMs) based on the connection network \citep{Abara1989} and the time-space network \citep{Hane1995}, respectively. With the increasingly fierce market competition, practitioners and researchers are interested in enhancing the basic FAMs by incorporating real-life considerations \citep{Sherali2006}. \cite{Farkas1996} is one of the pioneering studies to incorporate itinerary-based demand in the time-space FAM to capture the effects of flight interdependencies (network effects). \cite{Ioachim1999} address a weekly FAM with schedule synchronization constraints. \cite{Barnhart2002} propose an itinerary-based fleet assignment model (IFAM), which captures network effects and estimates passenger spill and recapture more accurately. \cite{Lohatepanont2004} address integrated models and algorithms for schedule design and fleet assignment, where the effect of flight frequencies and departure times on spilled passenger transfer preferences are considered. \cite{Belanger2006} incorporate the fleet homogeneity into the weekly FAM. To consider demand uncertainty, \cite{Sherali2008} develop a two-stage stochastic mixed-integer programming model. The first stage only focuses on the family-level assignments, while detailed assignments are performed based on demand realization in the second stage. \cite{Barnhart2009} propose a subnetwork fleet assignment model (SFAM) to better approximate the revenue side of the objective function, where the model simultaneously assigns fleet types to one or more flights in each subnetwork. \cite{Sherali2010} present a model that integrates the schedule design and fleet assignment, which directly incorporates itinerary-based demand for multiple fare classes. \cite{Yan2022} study an integrated schedule design and fleet assignment problem with a passenger choice model for fare product selections.

\subsection{Crew pairing problem} \label{review_cpp}
Crew pairing problem (CPP) is the key step of the crew scheduling problem, and aims to find a minimum-cost set of pairings to cover each flight. The CPP is typically formulated as a set-partitioning problem (SPP) \citep{Chu1997,Barnhart2003,Kasirzadeh2015}. Other formulations are also proposed. \cite{Vance1997} address a duty-based formulation and break the decision process into two stages: cover flights with duties and then build pairings with selected duties. \cite{Desaulniers1997} design an integer non-linear multi-commodity network flow problem with additional resource variables, where the nonlinearities are isolated in the subproblems. \cite{Yen2006} propose a two-stage stochastic programming model to identify pairings that are less sensitive to schedule disturbances. \cite{Shebalov2006} put forward the concept of move-up crews, which are the crews that can potentially be swapped in operations. By increasing the number of move-up crews, the pairing solution will be more robust. \cite{AhmadBeygi2009} consider a new integer programming approach to generate pairings. \cite{Saddoune2013} introduce a rolling horizon approach to solve the monthly pairing problem. \cite{Haouari2019} develop a novel compact polynomial-sized non-linear formulation, which is solved using the reformulation linearization technique. \cite{Antunes2019} present a robust optimization model, where the delay propagation through crew connections is captured. \cite{Quesnel2020} address a crew pairing problem with complex features (CPPCF), which considers the crew preferences to create more suitable pairings.

\subsection{Integration of fleet assignment and crew pairing problem} \label{review_int}
The outcome of the fleet assignment stage decomposes the flight network into subnetworks according to fleet families, such that the crew scheduling problem is solved for each fleet family. This is partially reasonable as a single crew member is only qualified to operate aircraft in the same fleet family. However, solving the CSP after the FAP can be suboptimal because the fleet assignment solution can limit the set of feasible pairings from which we can choose. Thus, some studies dedicate to integrating the FAP and CSP (mainly CPP). \cite{Barnhart1998} develop a model that integrates the fleet assignment and an approximation of the crew pairing problem, which is based on the duty network. \cite{Klabjan2002,Sandhu2007} integrate the FAP and CPP with a plane-count constraint to ensure that the number of aircraft used by crews does not exceed the number of available aircraft. \cite{Gao2009} propose an integrated FAP and CPP model using crew connections instead of explicit pairings. Moreover, they limit the number of fleet types and crew bases in each airport to enhance the robustness of the schedule. Because the aircraft routing problem is not considered, these integrated models do not guarantee maintenance feasibility. To satisfy the maintenance requirements, few studies consider the full integration of FAP, ARP, and CPP \citep{Clarke1996,Papadakos2009,SalazarGonzalez2014,Cacchiani2017,Shao2017,Cacchiani2020}.

To summarize, the existing fleet assignment studies mainly focus on the daily and weekly schedules. The long-term schedule, which spans a season or the whole year, is often constructed by repeating the short-term schedule. Furthermore, the integrated optimization of crew and aircraft resources mainly considers the pairing-related constraints, where pairings generally span 1-5 days. To the best of our knowledge, no research has incorporated the FAP with the macroscopical crew flight time allocation problem.

\section{Problem definition and formulation}\label{model}
Before introducing the integrated model of the TFACPP, we first briefly review two classic mathematical models for the fleet assignment and crew scheduling problems. Both models are essential components of the integrated model. Then, we introduce the basic integrated model.

\subsection{Fleet assignment problem} \label{model_fap}
Given a flight schedule and a set of aircraft, the FAP determines the fleet type for each flight leg with the objective of maximizing total profit. Typically, the fleet assignment model (FAM) is constructed based on a flight network with flight cover, flow balance, and aircraft number constraints. Next, we first introduce the widely used time-space network. Then, we describe the classic model based on this network, which is first proposed by \cite{Hane1995}.

\begin{figure}[h]
	\centering		
	\includegraphics[scale=0.5]{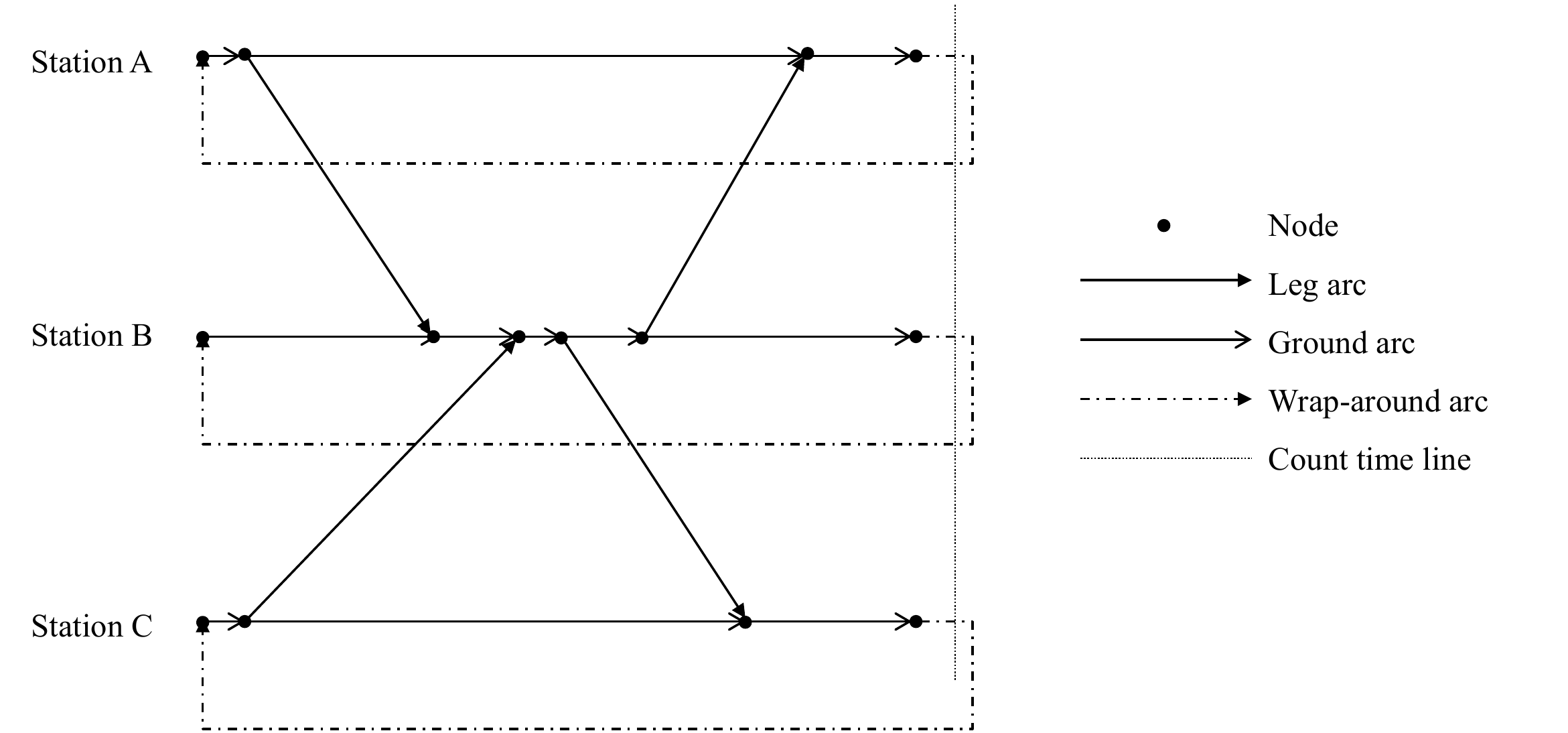}	
	\caption{\centering{An illustrate example of the time-space network}}		
	\label{fig_timespace}
\end{figure}

As shown in Fig. \ref{fig_timespace}, the time-space network is constructed for each fleet type and consists of a set of nodes representing departure or arrival events at stations, and a set of arcs representing flight legs, aircraft on the ground, and overnighting aircraft at stations. Specifically, the time of a departure event is the original departure time of the flight, whereas the time of an arrival event is the sum of the actual arrival time and minimum turn time for the aircraft. Finally, we arbitrarily select a point in time as the count time and use it to count the used aircraft number.

With the notations presented in Table \ref{tab_fampara}, the classic fleet assignment model is written as in Eqs. (\ref{fam_obj})–(\ref{fam_y}).

\begin{table}[h]
	\caption{Notations for the classic FAP}
	\centering
	\small
	{\def\arraystretch{1}  
		\begin{tabular}{lp{12.5cm}}
			\toprule
			\multicolumn{2}{l}{Sets and parameters} \\
			\midrule
			$L$	& Set of legs (leg arcs) in the flight schedule, indexed by $l$.\\
			$F$	& Set of fleet types in the airline, indexed by $f$, each fleet type is characterized by the number of aircraft $k_f$.\\
			${N_f}$ & Set of nodes in the time-space network for fleet type $f$, indexed by $n$.\\
			${G_f}$ & Set of ground arcs in the time-space network for fleet type $f$, indexed by $g$.\\
			$L_{nf}^ + \left( {L_{nf}^ - } \right)$ & Set of leg arcs beginning (ending) at node $n \in {N_f} \; \left( {\forall f \in F} \right)$.\\
			$G_{nf}^ + \left( {G_{nf}^ - } \right)$ & Set of ground arcs beginning (ending) at node $n \in {N_f} \; \left( {\forall f \in F} \right)$.\\
			$LC_f$ & Set of leg arcs passing the count timeline in the time-space network for fleet type $f$.\\
			$GC_f$ & Set of ground arcs passing the count timeline in the time-space network for fleet type $f$.\\
			${r_{lf}}$ & Profit of covering leg $l$ with aircraft of fleet type $f$, which is calculated by the difference between revenue and cost. \\
			\midrule
			Variables	&	\\
			\midrule
			${x_{lf}}$ & If leg $l$ is assigned to fleet type $f$, ${x_{lf}} = 1$; otherwise, ${x_{lf}} = 0$. \\
			${y_{gf}}$ & Number of aircraft of fleet type $f$ on the ground arc $g \in {G_f}$.\\
			\bottomrule
		\end{tabular}
		\label{tab_fampara}
	}
\end{table}

\begin{align}
	\text{Max} \hspace{0.25cm} & \sum\limits_{l \in L} {\sum\limits_{f \in F} {{r_{lf}}{x_{lf}}} } \label{fam_obj} \\
	\text{s.t.} \hspace{0.25cm} &  \sum\limits_{f \in F} {{x_{lf}}}  = 1 \;\;\;\; \forall l \in L \label{fam_cover}\\
	\ &  \sum\limits_{l \in L_{nf}^ + } {{x_{lf}}}  + \sum\limits_{g \in G_{nf}^ + } {{y_{gf}}}  = \sum\limits_{l \in L_{nf}^ - } {{x_{lf}}}  + \sum\limits_{g \in G_{nf}^ - } {{y_{gf}}} \;\;\;\; \forall f \in F,\forall n \in {N_f} \label{fam_blc} \\
	\ & \sum\limits_{l \in LC_f} {{x_{lf}}}  + \sum\limits_{g \in GC_f} {{y_{gf}}}  \le {k_f} \;\;\;\; \forall f \in F \label{fam_afn} \\
	\ & {x_{lf}} \in \left\{ {0,1} \right\} \;\;\;\; \forall l \in L,\forall f \in F \label{fam_x}\\
	\ & {y_{gf}} \ge 0 \;\;\;\; \forall f \in F,\forall g \in {G_f} \label{fam_y}
\end{align}

The objective function (\ref{fam_obj}) maximizes the total profit. Cover constraints (\ref{fam_cover}) ensure that each leg is covered by one and only one fleet type. Flow balance constraints (\ref{fam_blc}) ensure that the total inflow equals the total outflow for all nodes. Aircraft number constraints (\ref{fam_afn}) ensure that the number of used aircraft of each fleet type in the counting time does not exceed the available number. Constraints (\ref{fam_x}) and (\ref{fam_y}) are variable constraints. Especially, with the flow balance constraints (\ref{fam_blc}), integrality of $x$ variables implies integrality of $y$ variables.

\subsection{Crew pairing problem} \label{model_cpp}
The CPP is to find a minimum-cost set of crew pairings such that each flight leg is exactly covered by one pairing. A crew pairing is a legal sequence of flight legs that begins and ends at the same crew base. Generally, it can span 1-5 days, during which crew members are subject to various work-rule requirements established by the administration and union contracts. Because cockpit crews are qualified to serve only one fleet family (several fleet types require the same cockpit configuration and qualification), the CPP is usually solved separately for each fleet family.

With the additional notations presented in Table \ref{tab_cpppara}, the classic crew pairing model is written as in Eqs. (\ref{cpp_obj})–(\ref{cpp_z}).

\begin{table}[h]
	\caption{Additional notations for the classic CPP}
	\centering
	{\def\arraystretch{1}  
		\begin{tabular}{ll}
			\toprule
			\multicolumn{2}{l}{Sets and parameters} \\
			\midrule
			$P$	& Set of legal crew pairings, indexed by $p$.\\
			${c_p}$ & Cost of pairing $p$, which is mainly the salary of the crew members. \\
			${\delta _{lp}}$ & If leg $l$ is covered by pairing $p$, ${\delta _{lp}} = 1$; otherwise, ${\delta _{lp}} = 0$.	\\ 
			\midrule
			Variables	&	\\
			\midrule
			${z_p}$ & If pairing $p$ is included in the solution, ${z_p} = 1$; otherwise, ${z_p} = 0$.\\
			\bottomrule
		\end{tabular}
		\label{tab_cpppara}
	}
\end{table}

\begin{align}
	\text{Min} \hspace{0.25cm} & \sum\limits_{p \in P} {{c_p}{z_p}}  \label{cpp_obj} \\	
	\text{s.t.} \hspace{0.25cm} &  \sum\limits_{p \in P} {{\delta _{lp}}{z_p}}  = 1 \;\;\;\; \forall l \in L \label{cpp_cover}\\
	\ & {z_p} \in \left\{ {0,1} \right\} \;\;\;\; \forall p \in P \label{cpp_z}
\end{align}

The objective function (\ref{cpp_obj}) minimizes the cost of the chosen set of pairings. Constraints (\ref{cpp_cover}) and (\ref{cpp_z}) ensure that each leg is exactly covered by a pairing.

\subsection{The TFACPP formulation} \label{model_TFACPP}
The TFACPP aims to simultaneously assign a fleet type and a pairing to each flight, with the objective of maximizing the yearly profit minus the crew cost. The classic models and existing literature neither address the macroscopical crew consideration nor capture the demand fluctuation between months. In order to incorporate these concerns, we extend the single daily schedule to twelve flight schedules representing twelve months. With the additional notations in Table \ref{tab_basepara}, the complete model can be written as in Eqs. (\ref{merge_obj})–(\ref{merge_z}).

\begin{table}[h]
	\caption{Additional notations for the TFACPP formulation}
	\centering
	\small
	{\def\arraystretch{1}  
		\begin{tabular}{p{1.8cm}p{13.6cm}}
			\toprule
			\multicolumn{2}{l}{Sets and parameters} \\
			\midrule
			$ M $	&	Set of months, indexed by $ m $.\\
			$ L^m $	&	Set of legs (leg arcs) in month $ m $, indexed by $ l $. Here, we assume flights with the same flight number in different months are different flights. That is ${L^{{m_1}}} \cap {L^{{m_2}}} = \emptyset \; \left( {\forall {m_1},{m_2} \in M} \right)$. \\
			$ N_f^m $	&	Set of nodes in the time-space network for month $ m $ and fleet type $ f $, indexed by $ n $.	\\
			$ G_f^m $	&	Set of ground arcs in the time-space network for month $ m $ and fleet type $ f $, indexed by $ g $.	\\
			$L_{nf}^{m+} (L_{nf}^{m-} )$	&	Set of leg arcs beginning (ending) at node $ n \in N_f^m \; \left( {\forall m \in M,\; \forall f \in F} \right) $.	\\
			$G_{nf}^{m+} (G_{nf}^{m-} )$	&	Set of ground arcs beginning (ending) at node $ n \in N_f^m \; \left( {\forall m \in M,\; \forall f \in F} \right) $.	\\
			$LC_{f}^{m}$	&	Set of leg arcs passing the count timeline in the time-space network for month $ m $ and fleet type $ f $.	\\
			$GC_{f}^{m}$	&	Set of ground arcs passing the count timeline in the time-space network for month $ m $ and fleet type $ f $.	\\
			$ B $	&	Set of fleet families, indexed by $b$. Each fleet family contains $k_{b}$ numbers of crew members, which can operate one or more fleet types in the fleet family.	\\
			$ F_{b} $	& Set of fleet types in fleet family $b$. 	\\
			$ P_{b}^m $	& Set of pairings for fleet family $b$ and month $m$, indexed by $p$. Each set $ P_{b}^m $ is constructed using the complete set of $ L^m $. Specifically, we regard the same flight sequences for different fleet families as different pairings. \\
			$ r_{lf}^m $	&	Profit of covering leg $l \in {L_m}\;\left( {\forall m \in M} \right)$ with aircraft of fleet type $f$.	\\
			$ c_{pb}^m $	&	Cost of pairing $p \in P_b^m \; \left( {\forall b \in B,\forall m \in M} \right)$.	\\
			$t_{b} $ & Total available yearly crew flight time for fleet family $b$, which is calculated as ${k_{b}} \times \bar t _{b}$ (${k_{b}}$ is the crew number and $\bar t _{b}$ is the upper bound of yearly crew flight time per crew).\\
			$t_{b}^m $ & Total available monthly crew flight time for month $m$ and fleet family $b$, which is calculated as ${k_{b}} \times \bar t _{b}^m$ (${k_{b}}$ is the crew number and $\bar t _{b}^m$ is the upper bound of monthly crew flight time per crew). \\
			$t_{b}^l $ & Flight time of leg $l \in L_m \; \left( {\forall m \in M} \right)$ with fleet family $b$. The value is calculated by the product of frequency and flight duration, which depends partly on the fleet family. For example, it is possible that a flight operated by Boeing 777 is faster than by Boeing 737.\\
			$ t_{b}^p $	&	Flight time of pairing $p \in P_{b}^m \; \left( {\forall m \in M} \right)$, and $t_{b}^p = \sum\nolimits_{l \in L^m} {\delta_{lp} t_{b}^l} $.	\\	
			\midrule
			Variables	&	\\
			\midrule
			$ x_{lf}^m $	&	If leg $l \in {L_m}\;\left( {\forall m \in M} \right)$ is assigned to fleet type $ f $, ${x_{lf}^m} = 1$; otherwise, ${x_{lf}^m} = 0$.	\\
			$ y_{gf}^m $	&	Number of aircraft of fleet type $ f $ in the ground arc $g \in {G_{f}^m}\;\left( {\forall m \in M} \right)$.	\\
			$ {z_{pb}^m} $ & If pairing $p \in P_b^m \; \left( {\forall b \in B,\forall m \in M} \right)$ is included in the solution, $ {z_{pb}^m} = 1$; otherwise, $ {z_{pb}^m} = 0$. \\
			\bottomrule
		\end{tabular}
		\label{tab_basepara}
	}
\end{table}

\begin{align}
	 \text{Max} \hspace{0.25cm} & \sum\limits_{m \in M} {\sum\limits_{l \in {L^m}} {\sum\limits_{f \in F} {{r_{lf}^m}{x_{lf}^m}} } }  - \sum\limits_{m \in M} {\sum\limits_{b \in B} {\sum\limits_{p \in P_b^m} {{c_{pb}^m}{z_{pb}^m}} } }  \label{merge_obj} \\
	\text{s.t.} \hspace{0.25cm} & \sum \limits_{f \in F} {x_{lf}^m} = 1 \;\;\;\; \forall m \in M,\forall l \in {L^m}\label{basic_cover} \\
	\ & \sum \limits_{l \in L_{nf}^{m+} } {x_{lf}^m} + \mathop \sum \limits_{g \in G_{nf}^{m+} } {y_{gf}^m} = \sum \limits_{l \in L_{nf}^{m-} } {x_{lf}^m} + \sum \limits_{g \in G_{nf}^{m-} } {y_{gf}^m} \;\;\;\; \forall m \in M,\forall f \in F,\forall n \in {N_f^m}	\label{basic_blc} \\
	\ & \sum \limits_{l \in LC_{f}^{m}} {x_{lf}^m} + \sum \limits_{g \in GC_{f}^{m}} {y_{gf}^m} \le {k_f} \;\;\;\; \forall m \in M,\forall f \in F \label{basic_afn} \\
	\ & \sum\limits_{p \in P_b^m} {\delta_{lp} {z_{pb}^m}}  = \sum\limits_{f \in {F_b}} {{x_{lf}^m}} \;\;\;\; \forall m \in M,\forall l \in {L^m},\forall b \in B \label{merge_link}\\
	\ & \sum\limits_{p \in P_b^m} {{t_b^p}{z_{pb}^m}}  \le t_b^m \;\;\;\; \forall m \in M,\forall b \in B \label{merge_month} \\
	\ & \sum\limits_{m \in M} {\sum\limits_{p \in P_b^m} {{t_b^p}{z_{pb}^m}} }  \le t_b \;\;\;\; \forall b \in B \label{merge_year} \\
	\ & {x_{lf}^m} \in \left\{ {0,1} \right\}\;\;\;\;\forall m \in M,\forall l \in {L^m},\forall f \in F \label{basic_x} \\
	\ & {y_{gf}^m} \ge 0 \;\;\;\; \forall m \in M,\forall f \in F,\forall g \in {G_{f}^m} \label{basic_y} \\
	\ & {z_{pb}^m} \in \left\{ {0,1} \right\} \;\;\;\; \forall m \in M,\forall b \in B,\forall p \in P_b^m \label{merge_z}	
\end{align}

The objective function (\ref{merge_obj}) maximizes the difference between yearly profit and crew cost. Constraints (\ref{basic_cover})–(\ref{basic_afn}) are the counterpart of constraints (\ref{fam_cover})–(\ref{fam_afn}) in the classic FAP models. Constraints (\ref{merge_link}) link the FAP and CPP, ensuring that a pairing exactly covers a leg assigned to fleet types in the fleet family. Crew flight time constraints (\ref{merge_month}) and (\ref{merge_year}) ensure that the used crew flight time by all the selected pairings is less than the total available crew flight time for every month and the whole year. Constraints (\ref{basic_x})--(\ref{merge_z}) are the variable constraints. 

As shown in Eqs. (\ref{base_month})–(\ref{base_year}), because $t_{b}^p = \sum\nolimits_{l \in L^m} \delta_{lp} t_{b}^l \; \left( {\forall m \in M} \right) $ and the existing of linking constraints (\ref{merge_link}), we can substitute $\sum\nolimits_{p \in P_b^m} {{\delta _{lp}}{z_{pb}^m}}$ in the left hand side of monthly crew flight time constraints (\ref{merge_month}) with a leg-based form $\sum\nolimits_{l \in {L^m}} {\sum\nolimits_{f \in {F_b}} {t_b^l{x_{lf}^m}} }$. Similarly, we can replace $\sum\nolimits_{m \in M} {\sum\nolimits_{p \in P_b^m} {{t_b^p}{z_{pb}^m}} }$ in the left hand side of yearly crew flight time constraints (\ref{merge_year}) with the leg-based form $\sum\nolimits_{m \in M} {\sum\nolimits_{l \in {L^m}} {\sum\nolimits_{f \in {F_b}} {t_b^l{x_{lf}^m}} } }$.

\begin{align}
	\ & \sum\limits_{p \in P_b^m} {t_b^p{z_{pb}^m}} = \sum\limits_{p \in P_b^m}\sum\limits_{l \in L^m} {\delta_{lp}t_b^l{z_{pb}^m}} = \sum\limits_{l \in L^m} {t_b^l}\sum\limits_{p \in P_b^m}{\delta_{lp}z_{pb}^m} = \sum\limits_{l \in {L^m}} {\sum\limits_{f \in {F_b}} {t_b^l{x_{lf}^m}} }  \label{base_month}\\
	\ &  \sum\limits_{m \in M} {\sum\limits_{p \in P_b^m} {t_b^p{z_{pb}^m}} } = \sum\limits_{m \in M} {\sum\limits_{p \in P_b^m} {\sum\limits_{l \in L^m} {\delta_{lp}t_b^l{z_{pb}^m}}}} = \sum\limits_{m \in M}{\sum\limits_{l \in L^m} {t_b^l}} {\sum\limits_{p \in P_b^m} {\delta_{lp}{z_{pb}^m}} }  = \sum\limits_{m \in M} {\sum\limits_{l \in {L^m}} {\sum\limits_{f \in {F_b}} {t_b^l{x_{lf}^m}} } }  \label{base_year} 
\end{align}

Thus, constraints (\ref{merge_month})–(\ref{merge_year}) are equivalent to (\ref{basic_month})–(\ref{basic_year}), and the alternative basic integrated model (BIM) is written as follows.

\begin{align}
	\textbf{(BIM)} \;\;\;\; \text{Max} \hspace{0.25cm} & (\ref{merge_obj}) \nonumber \\
	\text{s.t.} \hspace{0.25cm} & (\ref{basic_cover}) – (\ref{merge_link}) \nonumber \\
	\ &  \sum \limits_{l \in {L^m}}  \sum \limits_{f \in {F_b}} {t_{b}^l}{x_{lf}^m} \le t_b^m \;\;\;\;\forall m \in M,\;\forall b \in B \label{basic_month} \\
	\ &  \sum \limits_{m \in M}  \sum \limits_{l \in {L^m}} \sum \limits_{f \in {F_b}} {t_{b}^l}{x_{lf}^m} \le t_b\;\;\;\;\forall b \in B \label{basic_year} \\
	\ & (\ref{basic_x}) – (\ref{merge_z}) \nonumber	
\end{align}

It is commonly believed that the CPP is much harder than the FAP. By replacing (\ref{merge_month})–(\ref{merge_year}) with (\ref{basic_month})–(\ref{basic_year}), we move the crew resource allocation decision from the crew pairing to the fleet assignment decision. So, we can decompose the model and decouple the CPP into different months.

\section{Solution method}\label{method}
The basic TFACPP model in Eqs. (\ref{merge_obj})–(\ref{merge_link}), (\ref{basic_month})–(\ref{basic_year}), and (\ref{basic_x})–(\ref{merge_z}) is too large to solve, where one of the difficulties lies in the exponential number of pairings and linking constraints (\ref{merge_link}). In practice, billions of possible pairings are common \citep{Barnhart1998}. To make matters worse, in the TFACPP, the crew schedule would work on the complete schedule instead of a separate schedule for each fleet family. Thus, we design a solution method based on the Benders decomposition, which results in a pair of problems that can be solved iteratively.


\subsection{The Benders reformulation}
For the given month $m$, fleet family $b$, and values $\bar x _{lf}^m$ (and $\bar y _{gf}^m$) satisfying fleet assignment constraints (\ref{basic_cover})--(\ref{basic_afn}) and macroscopical crew flight time constraints (\ref{basic_month})--(\ref{basic_year}), the linear programming (LP) relaxation of the basic formulation reduces to the following Benders primal subproblem (BSP):

\begin{align}
	\textbf{(BSP)} \;\;\;\; \text{Min} \hspace{0.25cm} &  {\sum\limits_{p \in P_b^m} {c_{pb}^mz_{pb}^m} } \label{BPSP_obj}  \\
	\text{s.t.} \hspace{0.25cm} & \sum\limits_{p \in P_b^m} {{\delta _{lp}}z_{pb}^m}  = \sum\limits_{f \in F_b} {\bar x_{lf}^m} \;\;\;\; \forall l \in {L^m} \label{BPSP_cover} \\
	\ & z_{pb}^m \ge 0 \;\;\;\; \forall p \in P_b^m \label{BPSP_z}
\end{align}

With cover constraints (\ref{basic_cover}), at most one of variables $\left\{ {x_{lf}^m:f \in {F_b}} \right\}$ is nonzero. More specifically, if a leg $l \in L^m$ is assigned to the fleet types in fleet family $b$, $\sum\nolimits_{f \in F_b} {\bar x_{lf}^m} = 1$; otherwise, $\sum\nolimits_{f \in F_b} {\bar x_{lf}^m} = 0$. Define $\bar L_b^m \subseteq L^m$ as the subset of legs in month $m$ that have been assigned to fleet types in fleet family $b$, constraints (\ref{BPSP_cover}) reduce to:

\begin{align}
 & \sum\limits_{p \in P_b^m} {{\delta _{lp}}z_{pb}^m}  = 1 \;\;\;\; \forall l \in {\bar L_b^m} \label{BPSP_simple_cover} 
\end{align}

The BSP in Eqs. (\ref{BPSP_obj}), (\ref{BPSP_simple_cover}), and (\ref{BPSP_z}) is the LP relaxation of the standard crew pairing problem, where the upper bounds on the $z_{pb}^m$ variables are not needed because of constraints (\ref{BPSP_simple_cover}) and the binary constraints on variables $z_{pb}^m$ have been relaxed to nonnegativity constraints. Define free variables $\omega _{lb}^m$ as the dual variables associated with cover constraints (\ref{BPSP_simple_cover}). The dual subproblem is written as:

\begin{align}
	\text{Max} \hspace{0.25cm} &  {\sum\limits_{l \in {\bar L_b^m}} { {\sum\limits_{f \in {F_b}} {\omega _{lb}^m} } } }  \label{BDSP_obj}  \\
	\text{s.t.} \hspace{0.25cm} & \sum\limits_{l \in {L_b^m}} {{\delta _{lp}}\omega _{lb}^m} \le c_{pb}^m \;\;\;\; \forall p \in P_p^m \label{BDSP_cap}	
\end{align}

A feasible solution to the BSP is easy to obtain, as long as a sufficient number of pairings are included. Thus, the TFACPP is a problem with complete recourse, and no feasibility cut is needed. With an additional free variable ${\eta}_b^m$, an optimality cut is written as (\ref{BMP_aoptcut}).

\begin{align}
	{\eta}_b^m  \ge  {\sum\limits_{l \in {\bar L_b^m}} { {\sum\limits_{f \in {F_b}} \bar \omega _{lb}^m{x_{lf}^m }}}} \;\;\;\; \forall m \in M, \forall b \in B  \label{BMP_aoptcut} 
\end{align}

Here, $\left( {\bar \omega  _{lb}^m|l \in {\bar L_b^m}} \right)$ is an extreme point of the dual polyhedron ${\Delta}_b^m $, defined by (\ref{BDSP_cap}) for month $m$ and fleet family $b$. With an additional set ${S_{\Delta}^{mb} }$, defined as the set of extreme points of ${\Delta}_b^m $, the BIM can be reformulated as:

\begin{align}
	\textbf{(BMP)} \;\;\;\; \text{Max} \hspace{0.25cm} & \sum\limits_{m \in M} {\sum\limits_{l \in {L^m}} {\sum\limits_{f \in F} {r_{lf}^mx_{lf}^m} } }  - \sum\limits_{m \in M} {\sum\limits_{b \in B}} {\eta}_b^m   \label{BMP_obj}  \\	
	\text{s.t.} \hspace{0.25cm} & {\eta}_b^m  \ge  {\sum\limits_{l \in {\bar L_b^m}} { {\sum\limits_{f \in {F_b}} {\bar \omega  _{lb}^mx_{lf}^m} } } } \;\;\;\;  {\forall m \in M}, {\forall b \in B}, \forall \left( {\bar \omega _{lb}^m|l \in {\bar L_b^m}} \right) \in {S_{\Delta}^{mb} } \label{BMP_optcut} \\
	\ & (\ref{basic_cover})–(\ref{basic_afn}), \; (\ref{basic_month})–(\ref{basic_year}), \; (\ref{basic_x})–(\ref{basic_y}) \nonumber
\end{align}

The above model is also called the Benders master problem (BMP), where each ${\eta}_b^m$ is restricted to be larger than or equal to the optimal value of the related subproblem by an optimality cut in (\ref{BMP_optcut}). Compared with the basic model, the BMP contains fewer variables but more constraints, although most optimality constraints are inactive in the optimal solution. Hence, instead of explicitly enumerating these constraints, it is common to use an iterative algorithm to generate only a subset of cuts that are sufficient to identify the optimal solution. More specifically, we start with a relaxation of the BMP without Benders cut and iteratively update the relaxed BMP with Benders cuts generated by solving the BSPs. This process terminates until the difference between the upper bound and lower bound is small enough.

\subsection{Generating Benders cuts with empirical data}
Traditionally, we solve the BSPs to generate the Benders cuts. As mentioned above, the BSP for every month and fleet family is the LP relaxation of a standard CPP, which can be solved using an existing algorithm such as column generation. Nevertheless, the procedure is computationally expensive. Moreover, our test cases contain 12 months and 6 fleet families, which indicates that 72 CPPs are required in an iteration. Thus, we analyze the economic meaning of the dual variables ${\omega}_{lb}^{m}$ and propose an alternative method to generate the Benders cuts. 

In the BSP for month $m$ and fleet family $b$, each constraint in Eqs. (\ref{BPSP_simple_cover}) represents a limitation for covering a leg $l \in \bar L_b^m$ with a pairing. The objective function (\ref{BPSP_obj}) represents the total pairing cost. Then, the dual variable $\omega _{lb}^m$ is the marginal cost of covering an additional leg. In other words, $\omega _{lb}^m$ is the fair price airline would pay to cover the additional leg with an arbitrary pairing. Then, the optimality cuts in Eqs. (\ref{BMP_optcut}) are to select a set of prices that maximizes its return. However, although the dual variables $\omega _{lb}^m$ provide accurate fair prices in mathematics, they underestimate the actual crew prices in fact. For example, the rostering-related cost is omitted in the TFACPP, where the roster is a monthly schedule for an individual crew member. Thus, we recommend replacing the dual variables $\omega _{lb}^m$ with given empirical data $\hat \omega _{lb}^m$. Then, we can substitute the accurate optimality cuts in Eqs. (\ref{BMP_optcut}) with approximate optimality cuts in Eqs. (\ref{BMP_approxcut}).

\begin{align}
	 \eta _b^m \ge \sum\limits_{l \in \bar L_b^m} {\sum\limits_{f \in {F_b}} {\hat \omega _{lb}^mx_{lf}^m} } \;\;\;\; \forall m \in M,\forall b \in B,\forall \bar L_b^m \subseteq {L^m} \label{BMP_approxcut}
\end{align}

The empirical value $\hat \omega _{lb}^m$ can be estimated by using existing/historical pairing and rosters to solve the crew scheduling problem. The number of existing/historical pairings and rosters is much less than the complete set of possible pairings and rosters, which are often of high quality. Thus, the dual information is easy to obtain.


\subsection{Solving the BMP with column generation}
As a significant number of pairing-related variables and constraints are isolated in the BSPs, the BMP is relatively solvable. However, there are still numerous leg-related variables and constraints. For industrial instances in our paper, the BMP comprises nearly 500,000 binary variables, which is almost impossible for a commercial solver to handle directly. Considering the special structure of the BMP, where the yearly constraints (\ref{basic_year}) serve as linking constraints and other constraints can be decomposed for each month, we apply a column generation algorithm to solve the problem efficiently.

\subsubsection{Column generation master problem} \label{method_mp}
Let ${{\rm{\Psi }}_m}$ be the set of all possible fleet assignment schedules for month $m \in M$, which meet all the monthly constraints. Each $\psi  \in {\Psi _m}$ is characterized by its profit $r_m^\psi  = \mathop \sum \limits_{l \in {L_m}} \mathop \sum \limits_{f \in F} {r_{lf}^m}x_{lf}^{m\psi} - \sum\limits_{b \in B} {\eta}_b^m $ and used crew flight time $t_{mb}^\psi  = \mathop \sum \limits_{l \in {L_m}} \mathop \sum \limits_{f \in {F_b}} {t_{b}^l}x_{lf}^{m\psi} \; \left( {\forall b \in B} \right)$. Moreover, we define binary variables $u_m^\psi $ to indicate whether the assignment schedule $\psi  \in {\Psi _m}$ is selected. Subsequently, we reformulate the BMP as the following column generation master problem (CGMP).

\begin{align}
	\textbf{CGMP} \;\;\;\; \text{Max} \hspace{0.25cm} & \sum \limits_{m \in M} \sum \limits_{\psi  \in {\Psi _m}} r_m^\psi u_m^\psi \label{mp_obj} \\
	\text{s.t.} \hspace{0.25cm} & \sum \limits_{\psi  \in {\Psi _m}} u_m^\psi  = 1\;\;\;\;\forall m \in M \label{mp_cover} \\
	\ & \sum \limits_{m \in M} \sum \limits_{\psi  \in {\Psi _m}} t_{mb}^\psi u_m^\psi  \le t_b\;\;\;\;\forall b \in B \label{mp_year} \\
	\ & u_m^\psi  \in \left\{ {0,1} \right\}{\rm{\;\;\;\;}}\forall m \in M,\;\forall \psi  \in {{\rm{\Psi }}_m} \label{mp_z}
\end{align}

The objective function (\ref{mp_obj}) aims to maximize the yearly profit, as in the basic formulation. Cover constraints (\ref{mp_cover}) ensure that only one fleet assignment schedule is selected for each month. Crew flight time constraints (\ref{mp_year}) ensure that the final schedule satisfies the yearly crew flight time limits.

\subsubsection{Column generation pricing subproblem} \label{method_sp}
We can obtain the optimal solution to the CGMP with all sets of monthly fleet assignment schedules. However, this solution method is impractical owing to the large number of ${{\rm{\Psi }}_m}$. Instead, we first initialize the CGMP with only a small subset of the monthly fleet assignment schedules, and then we iteratively update the CGMP by adding $u_m^\psi $ variables with positive reduced costs until we obtain the optimal LP solution.

Define ${\alpha _m}$ as the dual variables associated with cover constraints (\ref{mp_cover}), ${\beta _b}$ as the dual variables associated with yearly crew flight time constraints (\ref{mp_year}). The reduced cost of a monthly fleet assignment schedule, $\bar c_m^\psi $, is calculated as follows:

\begin{equation}
	\bar c_m^\psi  = r_m^\psi  - {\alpha _m} - \sum \limits_{b \in B} t_{mb}^\psi {\beta _b} = \sum \limits_{l \in {L^m}}  \sum \limits_{b \in B} \sum \limits_{f \in {F_b}} \left( {{r_{lf}^m} - {t_{b}^l}{\beta _b}} \right)x_{lf}^{m\psi} - \sum\limits_{b \in B} {\eta}_b^m - {\alpha _m} \label{rdc_u}
\end{equation}

Here, ${\alpha _m}$ is a constant for the certain month. Therefore, the optimal monthly fleet assignment schedule can be found through the following column generation pricing subproblem (CGSP).

\begin{align}
	\textbf{(CGSP)} \;\;\;\; \text{Max} \hspace{0.25cm} & {\chi _m} = \sum \limits_{l \in {L_m}} \sum \limits_{b \in B} \sum \limits_{f \in {F_b}} \left( {{r_{lf}^m} - {t_{b}^l}{\beta _b}} \right){x_{lf}^m} - \sum\limits_{b \in B} {{\eta}_b^m} \label{sp_obj} \\
	\text{s.t.} \hspace{0.25cm} & \sum \limits_{f \in F} {x_{lf}^m} \le 1\;\;\;\;\forall l \in {L_m} \label{sp_cover} \\
	\ & \sum \limits_{l \in L_{nf}^{m+} } {x_{lf}^m} + \sum \limits_{g \in G_{nf}^{m+} } {y_{gf}^m} = \sum \limits_{l \in L_{nf}^{m-} } {x_{lf}^m} + \sum \limits_{g \in G_{nf}^{m-} } {y_{gf}^m}\;\;\;\;\forall f \in F,\;\forall n \in {N_{f}^m} \label{sp_blc} \\
	\ & \sum \limits_{l \in L_{f}^{mk}} {x_{lf}^m} + \sum \limits_{g \in G_{f}^{mk}} {y_{gf}^m} \le {k_f}\;\;\;\;\forall f \in F \label{sp_afn} \\
	\ & \sum \limits_{l \in {L^m}} \sum \limits_{f \in {F_b}} {t_{b}^l}{x_{lf}^m} \le t_b^m\;\;\;\;\forall b \in B \label{sp_month} \\
	\ & {\eta}_b^m  \ge  {\sum\limits_{l \in {\bar L_b^m}} { {\sum\limits_{f \in {F_b}} {\hat \omega  _{lb}^mx_{lf}^m} } } } \;\;\;\; {\forall b \in B}, {\forall \bar L_b^m \subseteq {L^m}}  \label{sp_optcut} \\
	\ & {x_{lf}^m} \in \left\{ {0,1} \right\}\;\;\;\;\forall l \in {L^m},\;\forall f \in F \label{sp_x} \\
	\ & {y_{gf}^m} \ge 0\;\;\;\;\forall f \in F,\;\forall g \in {G_{f}^m} \label{sp_y} 
\end{align}

The objective function (\ref{sp_obj}) aims to maximize the reduced cost, while constraints (\ref{sp_cover})–(\ref{sp_y}) are analogous to monthly constraints in the BMP formulation. If the objective value ${\chi _m} > {\alpha _m}$, a new monthly fleet assignment schedule is found, and the corresponding variable  is added to the CGMP. 

The structure of the CGSP is similar to the FAM, with three differences. First, two additional parts are included in the objective function, which serve as the penalty of using the yearly crew flight time ($t_b^l \beta_b x_{lf}^{m}$) and the monthly crew cost ($\eta_b^m$). Second, we limit the monthly crew flight time with constraints (\ref{sp_month}). Third, we add the Benders cuts (\ref{sp_optcut}) to calculate the crew cost. In summary, the CGSP can be viewed as the decoupled monthly subproblem of the BMP, where the yearly crew flight time constraints are dualized to the objective function.

\subsubsection{Solution framework for CGMP}
Solving CGSP as an MIP problem could be very time-consuming. In our preliminary studies, it could take nearly one hour to obtain a good MIP solution. Fortunately, it is widely reported that the LP relaxation of the FAM-based model provides a good approximation of the MIP solution, because the majority of constraints are the network constraints \citep{Charnes1993,Hane1995,Lima2022}. 

\begin{figure}[h]
	\centering		
	\includegraphics[scale=0.65]{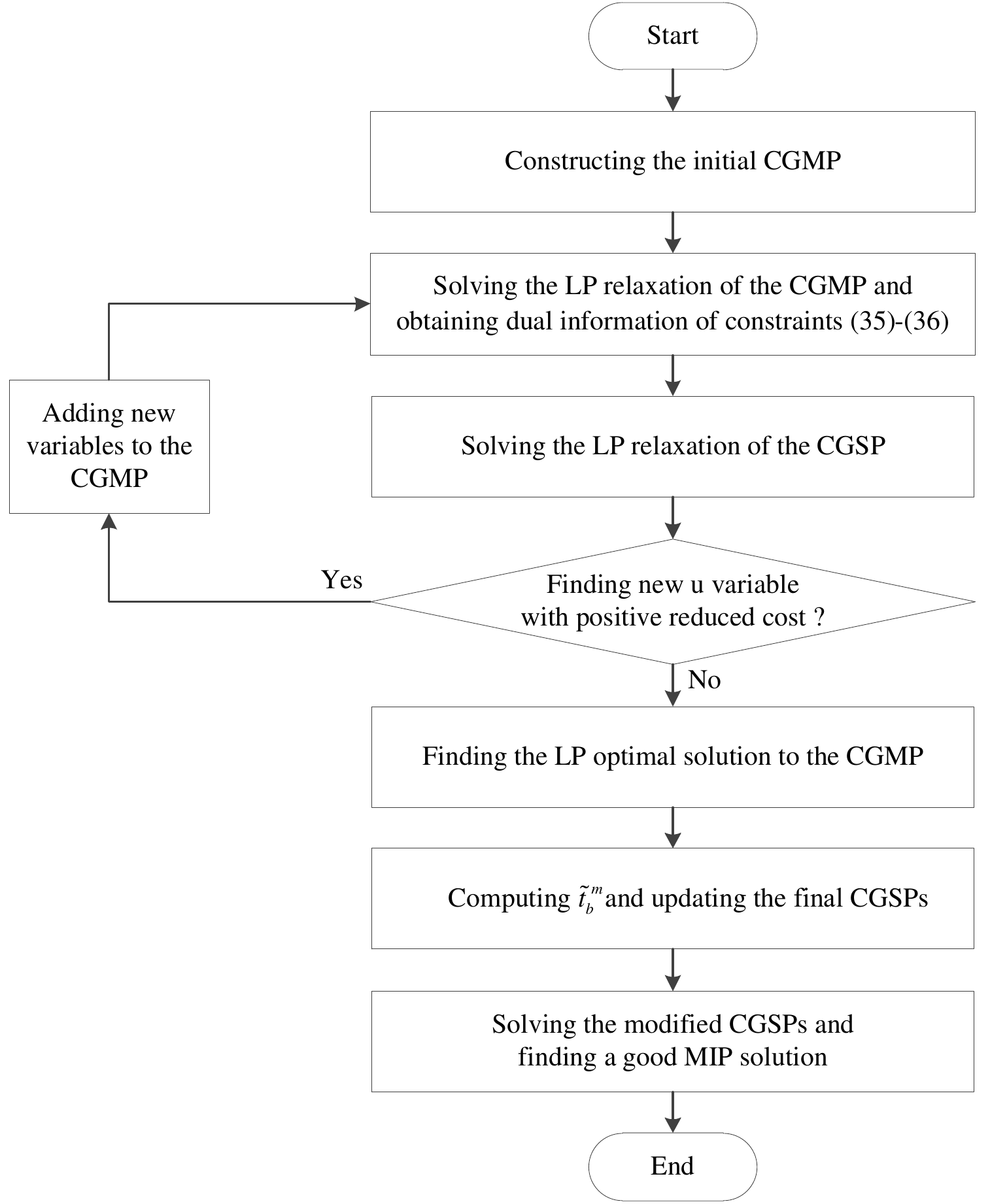}	
	\caption{\centering{Flow chart of the solution method to CGMP}}	
	\label{fig_flowchart}
\end{figure}

Consequently, we first obtain the optimal solution to the LP relaxation of the CGMP through column generation. Then, we solve the final CGSPs as MIP problems to obtain the binary fleet assignment decisions. It is important to note that if we solve the CGSPs directly, the accumulated yearly crew flight time of any fleet family may exceed the yearly crew flight time limitations. To resolve the problem, we modify the right hand side values of constraints (\ref{sp_month}) to be $\tilde t_b^m$, where $\tilde t_b^m = \sum\nolimits_{\psi  \in {\Psi _m}} {t_{mb}^\psi u_{mb}^\psi } $ in the LP optimal solution to the CGMP. It is easy to see that, $\tilde t_b^m$ is the optimal allocation of the yearly crew flight time for fleet family $b$ in month $m$ in the LP solution. The flow chart of the solution method is shown in Fig. \ref{fig_flowchart}.

\section{Extensions}\label{extend}
Airline daily operations and tactical management involve several practical considerations and additional requirements. As described in this section, we take into account two types of real considerations in the TFACPP formulation, so as to enhance the practicability of our problem. We first introduce the crew transition strategy and relative extensions. The strategy is well-practiced in the industry but seldom studied in the literature. Then, we incorporate the crew availability uncertainty into the TFACPP, which is caused by unexpected personnel absences such as illness and vacation.

\subsection{Extension to crew transition}\label{extend_tsf}
As mentioned previously, pilot development is very time-consuming and expensive. A typical trainee generally takes four to seven years to become a qualified captain. During this long duration, it is quite often for pilots to transfer from one fleet family to another \citep{Qi2004}. For example, pilots (even captains) for narrow-body fleet types are often chosen and trained for wide-body fleet types. Occasionally, when a new fleet type is acquired, airlines tend to quickly develop the pilot team through transition for safety reason and economic considerations. After all, the top priority of the airline is to smoothly operate the new aircraft as soon as possible. Similarly, when an old fleet type is retired, airlines also tend to relocate these pilots through transition. To summarize, the limited crew resource would be used more efficiently through crew transition. However, despite the importance of transition, decisions such as how many pilots to transfer are not well studied in the literature, and the relative decisions are often made manually in practice \citep{Akyurt_2021}. Thus, we extend the TFACPP with crew transition (TFACPP-CT) to provide better reference for transition decisions. The additional notations are listed in Table \ref{tab_tsfparas}. Using these notations, the TFATAP-CT can be formulated as the following MIP model in (\ref{tsf_obj})–(\ref{tsf_k}).

\begin{table}[h]
	\caption{Additional notations for the TFACPP-CT}
	\centering
	{\def\arraystretch{1}  
		\begin{tabular}{lp{11cm}}
			\toprule
			Parameters	&	\\
			\midrule
			$c_{{b_1}}^{{b_2}}$	&	Cost of transferring a crew member from fleet family $ b_1 $ to $ b_2 $.	\\
			$k_{{b_1}}^{{b_2}}$	&	Maximum number of crew members can be transferred from fleet family $ b_1 $ to $ b_2 $.	\\
			\midrule
			Variables	&	\\
			\midrule
			$v_{{b_1}}^{{b_2}}$	&	Number of crew members that plan to be transferred from fleet family $ b_1 $ to $ b_2 $.	\\
			\bottomrule
		\end{tabular}
		\label{tab_tsfparas}
	}
\end{table}

\begin{align}
	\textbf{(TFACPP-CT)} \;\;\;\; \text{Max} \hspace{0.25cm} & \sum \limits_{m \in M} \sum \limits_{l \in {L^m}} \sum \limits_{f \in F} {r_{lf}^m}{x_{lf}^m} -\sum\limits_{m \in M} \sum\limits_{b \in B} \sum\limits_{p \in P_b^m} c_{pb}^m z_{pb}^m - \sum \limits_{{b_1} \in B} \sum \limits_{{b_2} \in B\backslash {b_2}} c_{{b_1}}^{{b_2}}v_{{b_1}}^{{b_2}} \label{tsf_obj} \\
	\text{s.t.} \hspace{0.25cm} & (\ref{basic_cover}) – (\ref{merge_link}) \nonumber \\
	\ & \sum \limits_{l \in {L^m}} \sum \limits_{f \in {F_b}} {t_{b}^l}{x_{lf}^m} \le \left( {{k_b} + \sum \limits_{b' \in B\backslash b} v_{b'}^b - \sum \limits_{b' \in B\backslash b} v_b^{{b'}}} \right)\bar t_b^m\;\;\;\;\forall m \in M,\;\forall b \in B \label{tsf_month} \\
	\ & \sum \limits_{m \in M} \sum \limits_{l \in {L^m}} \sum \limits_{f \in {F_b}} {t_{b}^l}{x_{lf}^m} \le \left( {{k_b} + \sum \limits_{b' \in B\backslash b} v_{b'}^b - \sum \limits_{b' \in B\backslash b} v_b^{{b'}}} \right)\bar t_b\;\;\;\;\forall b \in B \label{tsf_year} \\
	\ & (\ref{basic_x}) – (\ref{merge_z}) \nonumber \\
	\ & v_{{b_1}}^{{b_2}} \in \left[ {0, k_{{b_1}}^{{b_2}}} \right]{\rm{\;\;\;\;}}\forall {b_1} \in B,{b_2} \in B\backslash {b_1} \label{tsf_k}
\end{align}

Compared with the original formulation, the TFATAP-CT has a different objective function and different available crew flight time. The objective function (\ref{tsf_obj}) aims to maximize the difference between the yearly profit and transition cost. The transition cost is mainly composed of the training cost and absence cost \citep{Qi2004}. The training cost is associated with the training processes, which are mandated by the governing administrations and labor unions. The absence cost is mainly the decrease in the profit owing to the reduced crew flight time. Once crews start the training process, they are unavailable for work until their training is complete. Generally, the absence cost is considerably large and hard to evaluate, which is an essential difficulty in crew transition planning. We estimate the absence cost by multiplying the training time and the yearly crew shadow price, which will be further introduced in Section \ref{insight}.

\subsection{Extension to crew uncertainty} \label{extend_rdm}

In the TFACPP, the available crew flight time is assumed to be deterministic. However, due to the uncertainty associated with illness, vacation, and training, crew availability is highly stochastic in practice \citep{Bayliss2017,Bayliss2020}, which will result in a lower utilization rate of the crew resource. The average yearly crew flight time for large-scale airlines could be less than 90\% of the upper limits. To ensure the feasibility of the operational stage, airlines need to compromise the crew resource in the planning stage.

To mitigate the impact of crew availability randomness, we extend the initial TFACPP by considering the crew flight time uncertainty (TFACPP-CU). Define ${\varepsilon _b} \in \left( {0,1} \right)$ as the risk tolerance of overusing yearly crew flight time, and $t_b\left( \rho  \right)$ as the stochastic yearly available crew flight time of fleet family $ b $. Then, the TFATAP-CU can be formulated as the following chance-constrained programming model.

\begin{align}
	\textbf{(TFACPP-CU)} \;\;\;\; \text{Max} \hspace{0.25cm} & \sum\limits_{m \in M} {\sum\limits_{l \in {L^m}} {\sum\limits_{f \in F} {r_{lf}^mx_{lf}^m} } } -\sum\limits_{m \in M} \sum\limits_{b \in B} \sum\limits_{p \in P_b^m} c_{pb}^m z_{pb}^m  \label{basic_obj}  \\
	\text{s.t.} \hspace{0.25cm} & (\ref{basic_cover}) – (\ref{merge_link}), \; (\ref{basic_month}) \nonumber \\
	\ & {\rm P}\left\{ { \sum \limits_{m \in M} \sum \limits_{l \in {L^m}} \sum \limits_{f \in {F_b}} {t_{b}^l}{x_{lf}^m} \le t_b\left( \rho  \right)} \right\} \ge 1 - {\varepsilon _b}\;\;\;\;\forall b \in B \label{rdm_year} \\
	\ & (\ref{basic_x}) – (\ref{merge_z}) \nonumber	
\end{align}

Constraints (\ref{rdm_year}) are individual chance constraints ensuring that the yearly crew flight time limit is satisfied with a probability of at least $1 - {\varepsilon _b}$ for each fleet family. The TFATAP-CU is a chance-constrained programming problem, which has been applied in various scenarios \citep{Shapiro2014,Yan2019}.

Without loss of generality, we assume that $t_b\left( \rho  \right)$ follows a finite discrete distribution for each fleet family. Specifically, we define ${Q_b}: = \left\{ {\rho _b^1, \ldots ,\rho _b^{{{\bar q}_b}}} \right\}$ as the finite ascending ordered set, the components of which are the possible realizations of the random variable $t_b\left( \rho  \right)$ with possibility $\phi _b^q = {\rm P}\left\{ {t_b\left( \rho  \right) = \rho _b^q} \right\}$ and $\mathop \sum \limits_{q \in {Q_b}} \phi _b^q = 1$. In this case, there exists an index ${q_0}$ such that constraint (\ref{rdm_q0}) is satisfied.

\begin{equation}
	 \sum \limits_{q = 1}^{{q_0} - 1} \phi _b^q < {\varepsilon _b} \le \mathop \sum \limits_{q = 1}^{{q_0}} \phi _b^q \label{rdm_q0}	
\end{equation}

The corresponding realization $\rho _b^{{q_0}}$ can be regarded as the $\left( {1 - {\varepsilon _b}} \right)$-quantile of the distribution $t_b\left( \rho  \right)$. Then, the deterministic equivalent of the individual chance constraints (\ref{rdm_year}) is as follows:

\begin{equation}
	 \sum \limits_{m \in M} \sum \limits_{l \in {L^m}} \sum \limits_{f \in {F_b}} {t_{b}^l}{x_{lf}^m} \le \rho _b^{{q_0}} \;\;\;\; \forall b \in B \label{rdm_demyear}	
\end{equation}

After conversion, the TFATAP-CU has the same structure with the basic TFATAP, but the usage of crew resource is more reasonable.

\section{Computational study}\label{test}
We first present the characteristics of the test cases considered in this study. Next, we describe the computational experiments of the TFACPP to which the Benders decomposition and column generation are applied. To highlight the benefit of the TFACPP, we compare the solution processes and results of the TFACPP with an equal allocation method. All the experiments are performed on a computer equipped with a 2.40 GHz CPU, 128 GB RAM, and 32 cores. All the mathematical models and algorithms are implemented using C++. Each LP problem is solved through a callable CPLEX library version 12.7 with the default settings.

\subsection{Test cases and parameter sets}\label{caseintro}
We assess our TFACPP formulation and solution method based on a yearly network, which is provided by a major Chinese airline. This network is the only large-scale network for which extensive demand data and resource data needed to solve the TFACPP (flight sets, fleet types, crew flight time, etc.) are available. To demonstrate the efficiency of the proposed model and solution method, we perturb the original demand to obtain two additional sets of demand for each flight leg. First, the demand for each leg is multiplied by a random parameter in the range [1.1,1.2], and the resulting set is defined as the high-demand level set. Next, the demand for each leg is multiplied by a random parameter in the range [0.8,0.9], and the resulting set is defined as the low-demand level set. Finally, the original set is defined as the mid-demand level set.

The whole year is divided into three time periods: $ m $ = 1–3, 4–10, and 11–12. The division is in line with the actual situation for most airlines: the winter and summer schedules are considerably different in terms of the number of flights and destinations. The summer ($ m $ = 4–10) flight set involves 1488 legs and 188 stations, whereas the winter ($ m $ = 1–3 and 11–12) flight set involves 1572 legs and 191 stations. Each flight leg has a specific origin, destination, departure time, arrival time, operation cost for every fleet type, the average demand for every month, average ticket price for every month, and flight time for every month and fleet family. These values are estimated from historical data using the airline’s forecasting systems.

The fleet consists of 380 airplanes, which are distributed in 22 fleet types and 6 fleet families. The seat number $ s_f $ and aircraft number $ k_f $ for each fleet type are specified in Table \ref{tab_fleet}. In particular, the fleet types A319-1 and A319-2 (B738-1 and B738-2) are similar except for the engine performance. Moreover, because pilots are usually the scarcest crew resources for airlines, we only consider the pilot flight time limitations in the computational experiments. The airline has more than 1700 pilots, among which the pilots for wide-body (A330, B747, B777, and B787) and narrow-body fleet types (A320 and B737) account for approximately 40\% and 60\% of the total, respectively.

\begin{table}[h]
	\caption{Information regarding fleet types}
	\centering
	{\def\arraystretch{1}  
		\begin{tabular}{p{2.15cm}p{2.15cm}rr}
			\toprule
			Fleet type & Fleet family & $ s_f $  & $ k_f $ \\
			\midrule
			A319-1 & A320  & 128   & 12 \\
			A319-2 & A320  & 128   & 21 \\
			A320  & A320  & 158   & 42 \\
			A321-1 & A320  & 177   & 29 \\
			A321-2 & A320  & 185   & 27 \\
			A332-1 & A330  & 237   & 18 \\
			A332-2 & A330  & 265   & 10 \\
			A332-3 & A330  & 283   & 4 \\
			A333-1 & A330  & 311   & 6 \\
			A333-2 & A330  & 301   & 13 \\
			A333-3 & A330  & 301   & 5 \\
			B737  & B737  & 128   & 21 \\
			B738-1 & B737  & 167   & 4 \\
			B738-2 & B737  & 167   & 57 \\
			B738-3 & B737  & 159   & 54 \\
			B738-4 & B737  & 176   & 14 \\
			B744  & B747  & 344   & 2 \\
			B748  & B747  & 365   & 6 \\
			B772  & B777  & 310   & 5 \\
			B773  & B777  & 311   & 20 \\
			B77W  & B777  & 392   & 3 \\
			B789  & B787  & 293   & 7 \\
			\bottomrule
		\end{tabular}%
		\label{tab_fleet}
	}
\end{table}

\subsection{Computational experiments of the TFACPP}\label{casemp}
We first solve the initial BMP or CGMP and present the results in this section. The monthly and yearly upper bounds of the flight time per crew, $\bar t_b^m$ and $\bar t_b$, are set as 100 h and 1000 h, respectively. This section presents the comparison of the CPU time and computational effort required to solve the CGMP for test cases involving three demand levels. As mentioned previously, we first solve the LP relaxation of the CGMP using the column generation algorithm. The computational results for each demand level (Dmd level) are presented in Table \ref{tab_lpmp}, including the number of calls to the LP relaxation of the CGMP (No. of CGMP), the average computation time of the LP relaxation of the CGMP in seconds (Avg. CGMP time), the number of calls to the LP relaxation of the CGSP (No. of CGSP), the average computation time of the LP relaxation of the CGSP in seconds (Avg. CGSP time), the number of generated monthly schedules (No. of sche), the objective value of the LP in billion (Yr LP obj), the total computation time for solving the LP relaxation in hours (Yr. LP time). It is worth to note that we solve the CGSP for each month without dual information to construct the initial CGMP.

\begin{table}[htbp]
	\caption{Computational results of the LP relaxation of the CGMP with column generation}
	\centering
	{\def\arraystretch{1}  
		
		\begin{tabular}{llllllll}
			\toprule
			\multicolumn{1}{p{1cm}}{Dmd level} & \multicolumn{1}{p{1.1cm}}{No. of CGMP } & \multicolumn{1}{p{2.1cm}}{Avg. CGMP time (s)} & \multicolumn{1}{p{1.1cm}}{No. of CGSP} & \multicolumn{1}{p{2.1cm}}{Avg. CGSP time (s)} & \multicolumn{1}{p{1cm}}{No. of sche} & \multicolumn{1}{p{1.5cm}}{Yr LP obj (bn)} & \multicolumn{1}{p{1.5cm}}{Yr LP time (h)} \\
			\midrule
			\multicolumn{1}{l}{High} & \multicolumn{1}{r}{41} & \multicolumn{1}{r}{0.000732} & \multicolumn{1}{r}{504} & \multicolumn{1}{r}{17.652} & \multicolumn{1}{r}{476} & \multicolumn{1}{r}{46.155 } & \multicolumn{1}{r}{3.933} \\
			\multicolumn{1}{l}{Mid} & \multicolumn{1}{r}{37} & \multicolumn{1}{r}{0.000649} & \multicolumn{1}{r}{456} & \multicolumn{1}{r}{15.758} & \multicolumn{1}{r}{413} & \multicolumn{1}{r}{39.105 } & \multicolumn{1}{r}{3.323} \\
			\multicolumn{1}{l}{Low} & \multicolumn{1}{r}{30} & \multicolumn{1}{r}{0.001667} & \multicolumn{1}{r}{372} & \multicolumn{1}{r}{16.475} & \multicolumn{1}{r}{337} & \multicolumn{1}{r}{31.655 } & \multicolumn{1}{r}{2.802} \\
			\bottomrule
		\end{tabular}%
		\label{tab_lpmp}%
	}
\end{table}

After obtaining the optimal LP solutions, we solve the modified CGSPs as MIP problems to obtain a good MIP solution. We indicate in Table \ref{tab_iphigh}–\ref{tab_iplow}, the objective value of the LP in billion (Mon LP obj), the objective value of the MIP in billion (Mon MIP obj), the integrality gap (Mon Int gap) defined as $\left( {Mon \;{\rm{ }}MIP \;{\rm{ }}obj{\rm{ }} - {\rm{ }}Mon \;{\rm{ }}LP \;{\rm{ }}obj} \right)/Mon \;{\rm{ }}LP \;{\rm{ }}obj * 100\%$, and the computation time of solving the MIP of the SPs in minutes (Mon MIP time) for each demand level and month.

\begin{table}[htbp]
	\caption{Computational results of the MIP CGSPs (High)}
	\centering
	{\def\arraystretch{1}  
		\begin{tabular}{rrrrr}
			\toprule
			\multicolumn{1}{p{1cm}}{Month} & \multicolumn{1}{p{1.5cm}}{Mon LP obj (bn)} & \multicolumn{1}{p{1.7cm}}{Mon MIP obj (bn)} & \multicolumn{1}{p{1.5cm}}{Mon Int gap} &  \multicolumn{1}{p{1.7cm}}{Mon MIP time (m)} \\
			\midrule
			1     & 3.549  & 3.547  & -0.057\% &  70.943 \\
			2     & 3.130  & 3.129  & -0.042\% &  193.138 \\
			3     & 2.869  & 2.868  & -0.038\% &  41.663 \\
			4     & 3.890  & 3.890  & -0.012\% &  51.100 \\
			5     & 3.577  & 3.576  & -0.030\% &  67.628 \\
			6     & 3.865  & 3.863  & -0.051\% &  63.553 \\
			7     & 5.232  & 5.230  & -0.051\% &  356.647 \\
			8     & 5.305  & 5.303  & -0.040\% &  324.885 \\
			9     & 4.267  & 4.265  & -0.046\% &  193.350 \\
			10    & 3.870  & 3.869  & -0.040\% &  177.893 \\
			11    & 3.293  & 3.292  & -0.029\% &  246.878 \\
			12    & 3.307  & 3.305  & -0.050\% &  276.237 \\
			The year & 46.155 & 46.137 & -0.041\% & 2063.915 \\
			\bottomrule
		\end{tabular}%
		\label{tab_iphigh}					
	}
\end{table}

\begin{table}[htbp]
	\caption{Computational results of the MIP CGSPs (Mid)}
	\centering
	{\def\arraystretch{1}  
		\begin{tabular}{rrrrr}
			\toprule
			\multicolumn{1}{p{1cm}}{Month} & \multicolumn{1}{p{1.5cm}}{Mon LP obj (bn)} & \multicolumn{1}{p{1.7cm}}{Mon MIP obj (bn)} & \multicolumn{1}{p{1.5cm}}{Mon Int gap} &  \multicolumn{1}{p{1.7cm}}{Mon MIP time (m)} \\
			\midrule
			1     & 2.933  & 2.930  & -0.090\% &  66.888 \\
			2     & 2.585  & 2.584  & -0.048\% &  45.235 \\
			3     & 2.386  & 2.386  & -0.035\% &  19.619 \\
			4     & 3.293  & 3.292  & -0.014\% &  42.362 \\
			5     & 2.996  & 2.996  & -0.025\% &  46.423 \\
			6     & 3.259  & 3.256  & -0.095\% &  81.549 \\
			7     & 4.509  & 4.506  & -0.059\% &  93.825 \\
			8     & 4.644  & 4.642  & -0.042\% &  104.665 \\
			9     & 3.650  & 3.648  & -0.042\% &  61.427 \\
			10    & 3.290  & 3.289  & -0.026\% &  10.264 \\
			11    & 2.760  & 2.759  & -0.034\% &  41.006 \\
			12    & 2.800  & 2.800  & -0.015\% &  24.939 \\	
			The year & 39.105 & 39.088 & -0.014\% & 638.202 \\		
			\bottomrule
		\end{tabular}%
		\label{tab_ipmid}					
	}
\end{table}

\begin{table}[htbp]
	\caption{Computational results of the MIP CGSPs (Low)}
	\centering
	{\def\arraystretch{1}  
		\begin{tabular}{rrrrr}
			\toprule
			\multicolumn{1}{p{1cm}}{Month} & \multicolumn{1}{p{1.5cm}}{Mon LP obj (bn)} & \multicolumn{1}{p{1.7cm}}{Mon MIP obj (bn)} & \multicolumn{1}{p{1.5cm}}{Mon Int gap} &  \multicolumn{1}{p{1.7cm}}{Mon MIP time (m)} \\
			\midrule
			1     & 2.338  & 2.337  & -0.059\% &  103.948 \\
			2     & 2.064  & 2.063  & -0.016\% &  47.550 \\
			3     & 1.929  & 1.929  & -0.022\% &  23.428 \\
			4     & 2.682  & 2.682  & -0.016\% &  25.152 \\
			5     & 2.392  & 2.391  & -0.064\% &  30.727 \\
			6     & 2.624  & 2.620  & -0.144\% &  27.256 \\
			7     & 3.663  & 3.661  & -0.070\% &  56.910 \\
			8     & 3.837  & 3.833  & -0.106\% &  53.447 \\
			9     & 2.966  & 2.962  & -0.165\% &  42.386 \\
			10    & 2.658  & 2.656  & -0.100\% &  32.856 \\
			11    & 2.231  & 2.230  & -0.028\% &  16.300 \\
			12    & 2.270  & 2.266  & -0.176\% &  30.954 \\	
			The year & 31.655 & 31.628 & -0.084\% & 490.913 \\		
			\bottomrule
		\end{tabular}%
		\label{tab_iplow}						
	}
\end{table}

For the yearly problem, the total computation time, including the solution of the LP relaxation of CGMP and MIP of CGSPs, is 38.333, 13.961, and 10.984 h in the high-, mid-, and low-demand levels, respectively. Because the TFACPP is a scheduling problem considered far in advance, the computation time and result quality are acceptable. Moreover, the computational results show that the integrality gaps between the MIP and LP objective values are nearly zero for all three demand levels. This finding shows that the proposed column generation solution method is efficient and effective for the CGMP.

However, the exhaustive solution method based on the Benders decomposition framework is very hard. Compared to the BSP, the expected computational efforts required to solve the BMP or CGMP will be much more. As shown in Table \ref{tab_lpmp}--\ref{tab_iplow}, the initial BMP is hard enough. With empirical data, the Benders cuts are easy to generate. Moreover, the Benders cuts do not affect the structure of the BMP. Thus, we only focus on the solution procedures for the initial BMP in this study. By preserving the macroscopical crew flight time limitations, we do not relax the important crew-related decisions.

\begin{figure}[htbp]
	\centering			
	\includegraphics[scale=0.7]{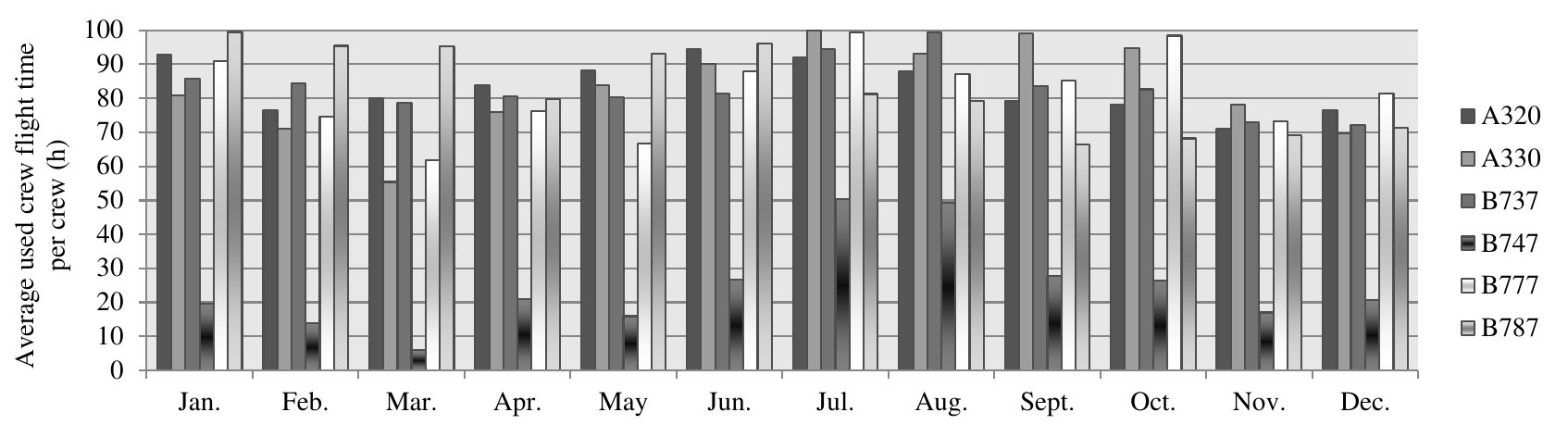}	
	\caption{\centering{Crew flight time allocation in the high-demand level}}
	\label{fig_hightime}
\end{figure}

\begin{figure}[htbp]
	\centering		
	\includegraphics[scale=0.7]{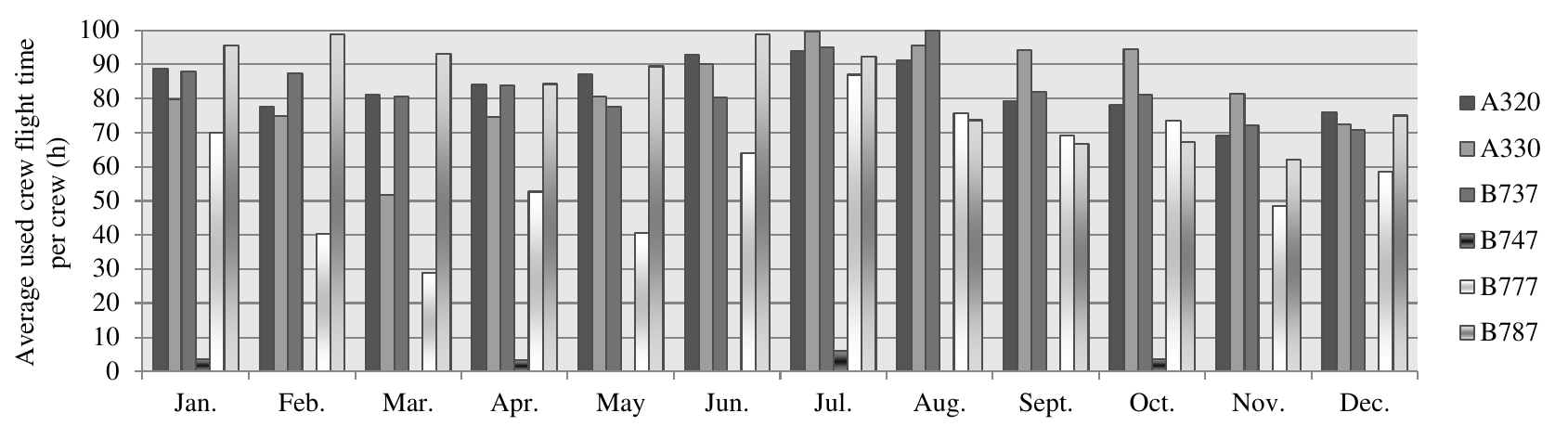}	
	\caption{\centering{Crew flight time allocation in the mid-demand level}}
	\label{fig_midtime}
\end{figure}

\begin{figure}[htbp]
	\centering			
	\includegraphics[scale=0.7]{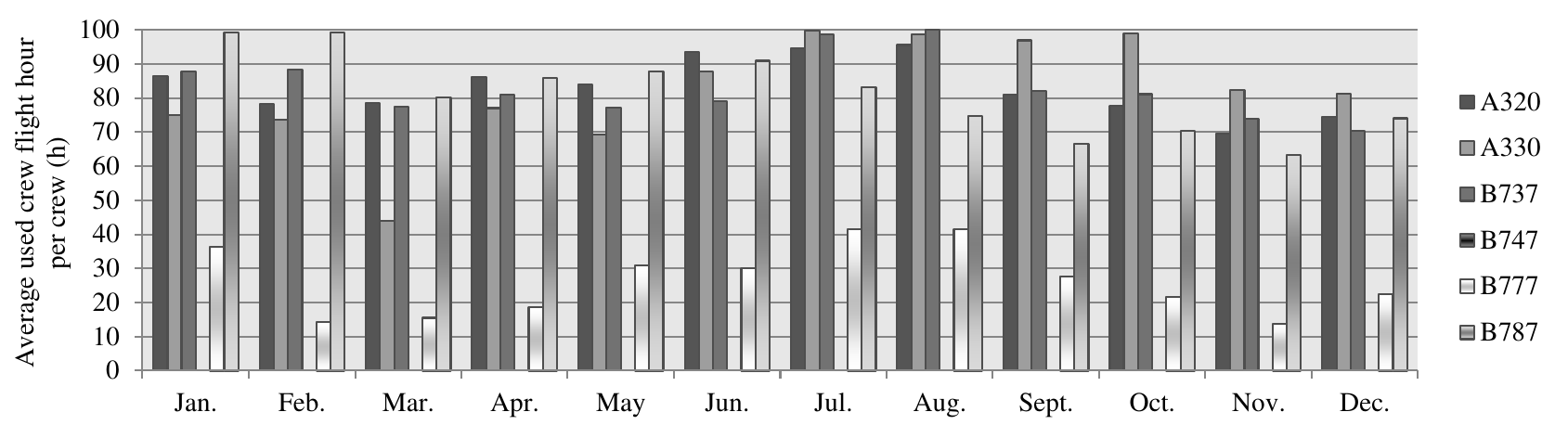}	
	\caption{\centering{Crew flight time allocation in the low-demand level}}
	\label{fig_lowtime}
\end{figure}

Furthermore, we investigate the used crew flight time in the approximate optimal MIP solutions. The considered details include the average used monthly crew flight time per crew for each fleet family as shown in Fig. \ref{fig_hightime}–\ref{fig_lowtime}, and the average used yearly crew flight time per crew for each fleet family as shown in Table \ref{tab_yeartime}. In Fig. \ref{fig_hightime}–\ref{fig_lowtime}, the horizontal and vertical coordinates indicate the month and average used monthly crew flight time per crew for each fleet family, respectively.

\begin{table}[htpb]
	\caption{Average used yearly crew flight time per crew (h)}
	\centering
	{\def\arraystretch{1}  
		\begin{tabular}{lrrrrrr}
			\toprule
			Dmd level & \multicolumn{1}{l}{A320} & \multicolumn{1}{l}{A330} & \multicolumn{1}{l}{B737} & \multicolumn{1}{l}{B747} & \multicolumn{1}{l}{B777} & \multicolumn{1}{l}{B787} \\
			\midrule
			High  & 997.347 & 995.522 & 996.360 & 292.714 & 983.906 & 973.923 \\
			Mid   & 996.399 & 993.093 & 998.149 & 16.638 & 711.666 & 974.452 \\
			Low   & 997.372 & 989.340 & 996.982 & 0.000 & 315.040 & 956.291 \\
			\bottomrule
		\end{tabular}%
		\label{tab_yeartime}					
	}
\end{table}

As shown in Fig. \ref{fig_hightime}–\ref{fig_lowtime} and Table \ref{tab_yeartime}, the monthly and yearly crew flight time limitations are satisfied in the computational results of the CGMP for all the demand levels. In addition, the allocation of the crew flight time resources is consistent with the airline trends, and is more scientific and accurate. For instance, for all three demand levels and fleet families, the maximum crew flight time resources are allocated to January–February and July–August, which are the traditional peak seasons. Furthermore, the following additional findings are derived, which may be helpful to the airline:

\begin{enumerate}
	\item The crew resources may be redundant for fleet families B747 and B777. For example, the crew flight time resources of B747 are surplus in all demand levels, and the crew flight time resources of B777 are surplus in the mid- and low-demand levels. The redundancy may have been observed by the airline, and the analysis here provides more scientific evidence.

	\item For most fleet families, the summer holiday season (July and August) is more important than the winter holiday season (January and February). Take B737 as an example, the differences in the average used crew flight time per crew between the summer and winter holiday seasons are 23.847, 19.826, and 22.384 h in the high-, mid-, and low-demand levels, respectively.

	\item The fleet family B787 is special, for which the winter holiday season (January and February) is more important than the summer holiday season (July and August). The differences in the average used crew flight time per crew between the winter and summer holiday seasons are 34.362, 28.227, and 40.574 h in the high-, mid-, and low-demand levels, respectively.
	
	\item The crew demand for narrow-body fleet families is more stable than that for wide-body fleet families. Consider, for example, the cases of B737 and A330 (for which the crew resources are almost completely used for the complete year): The average used monthly crew flight time of B737 is more than 70\% of the upper bound for all months; however, the corresponding value of A330 is less than 50\% of the upper bound for several months. Because large differences in crew resource usage adversely impact crew scheduling, measures should be adopted to boost the demand during the downturn.
	
	\item For fleet family A320, the demand in June should be seriously considered. For all demand levels, the crew flight time resources of A320 are almost exhausted, and the resources allocated to June rank at least third. Especially, June consumes the most resources among all months in the high-demand level.
\end{enumerate}

Other characteristics of the monthly used crew flight time are listed in Table \ref{tab_difftime}, which presents the maximum and minimum values of the average used monthly crew flight time per crew (Mon max, Mon min) and their differences (Diff).

\begin{table}[htpb]
	\caption{Difference in the average used monthly crew flight time (h)}
	\centering
	\small
	{\def\arraystretch{1}  
		\begin{tabular}{lrrrrrrrrr}
			\toprule
			\multicolumn{1}{c}{\multirow{2}[4]{*}{\begin{tabular}[c]{@{}l@{}}Dmd \\ level\end{tabular}}} & \multicolumn{3}{c}{High-demand level} & \multicolumn{3}{c}{Mid-demand level} & \multicolumn{3}{c}{Low-demand level} \\
			\cmidrule{2-10}          & \multicolumn{1}{l}{Mon max} & \multicolumn{1}{l}{Mon min} & \multicolumn{1}{r}{Diff} & \multicolumn{1}{l}{Mon max} & \multicolumn{1}{l}{Mon min} & \multicolumn{1}{r}{Diff} & \multicolumn{1}{l}{Mon max} & \multicolumn{1}{l}{Mon min} & \multicolumn{1}{r}{Diff} \\
			\midrule
			A320  & 94.454 & 70.932 & 23.522 & 93.990 & 69.100 & 24.889 & 95.549 & 69.560 & 25.989 \\
			A330  & 99.995 & 55.436 & 44.560 & 99.567 & 51.776 & 47.792 & 99.830 & 44.001 & 55.829 \\
			B737  & 99.406 & 72.244 & 27.162 & 99.969 & 70.851 & 29.118 & 99.966 & 70.478 & 29.487 \\
			B747  & 50.240 & 5.910 & 44.330 & 5.897 & 0.000 & 5.897 & 0.000 & 0.000 & 0.000 \\
			B777  & 99.204 & 61.685 & 37.519 & 86.969 & 28.752 & 58.217 & 41.594 & 13.771 & 27.823 \\
			B787  & 99.420 & 66.306 & 33.114 & 98.886 & 62.010 & 36.876 & 99.268 & 63.189 & 36.078 \\
			\bottomrule
		\end{tabular}%
		\label{tab_difftime}					
	}
\end{table}

The different distributions and large differences in the crew flight time resource usage demonstrate the validity of the proposed TFACPP.

\subsection{Comparison with an equal allocation method}\label{casecompare}
To examine the benefit of the proposed model, we compare the computational results of CGMP with an equal allocation method (EAM). In the EAM, we first equally allocate the yearly crew flight time into different months, and then solve the EAM CGSPs, where the right hand side value of constraints (\ref{sp_month}) are set to $\bar t_b^a = 1000/12 = 83.3$ h and $\beta_b = 0$. Since no dual information works as penalties of use the yearly crew flight time in the EAM CASPs, the monthly crew flight time constraints are hard to meet and the MIP solution is difficult to obtain. Thus, we only solve the LP relaxation here. Moreover, we view the values of $\sum\nolimits_{m \in M} \sum\nolimits_{l \in L^m} \sum\nolimits_{f \in F} r_{lf}^m x_{lf}^m$ in the LP optimal solution as the profit. The comparison of the computational results between these two models is presented in Table \ref{tab_compobj}, with the growth profit defined as (CGMP profit – EAM profit) and the growth rate defined as (CGMP profit – EAM profit)/EAM profit * 100\%. The proposed TFACPP formulation and solution method can increase the profit by more than two hundred million, which is a significant amount for the airline. The comparison of the average used yearly crew flight time per crew is presented in Table \ref{tab_comptime}. The differences in the average used yearly crew flight time between the two models are extremely small at all demand levels.

\begin{table}[htpb]
	\caption{Comparison of computation results}
	\centering
	{\def\arraystretch{1}  
		\begin{tabular}{lrrrrrr}
			\toprule
			\multicolumn{1}{c}{\multirow{2}[4]{*}{\begin{tabular}[c]{@{}l@{}}Dmd \\ level\end{tabular}}} & \multicolumn{2}{c}{EAM} & \multicolumn{2}{c}{CGMP} & \multicolumn{1}{c}{\multirow{2}[4]{*}{\begin{tabular}[c]{@{}l@{}}Growth \\ profit (bn)\end{tabular}}} &	
			 \multicolumn{1}{c}{\multirow{2}[4]{*}{\begin{tabular}[c]{@{}l@{}}Growth \\ rate \end{tabular}}}\\
			\cmidrule{2-5}          & \multicolumn{1}{l}{Profit (bn)} & \multicolumn{1}{l}{LP time (h)} & \multicolumn{1}{l}{Profit (bn)} & \multicolumn{1}{l}{LP time (h)} & & \\
			\midrule
			High  & 58.126 & 0.063  & 58.402 & 3.933 & 0.276 & 0.475\% \\
			Mid   & 49.348 & 0.059  & 49.564 & 3.323 & 0.216 & 0.438\% \\
			Low   & 38.656 & 0.063  & 38.876 & 2.802 & 0.220 & 0.570\% \\
			\bottomrule
		\end{tabular}%
		\label{tab_compobj}					
	}
\end{table}

\begin{table}[htpb]
	\caption{Comparison of average used yearly crew flight time per crew (h)}
	\centering
	\small
	{\def\arraystretch{1}  
		\begin{tabular}{lrrrrrrrrrrrr}
			\toprule
			\multicolumn{1}{c}{\multirow{2}[4]{*}{\begin{tabular}[c]{@{}l@{}}Dmd \\ level\end{tabular}}} & \multicolumn{2}{c}{A320} & \multicolumn{2}{c}{A330} & \multicolumn{2}{c}{B737} & \multicolumn{2}{c}{B747} & \multicolumn{2}{c}{B777} & \multicolumn{2}{c}{B787} \\
			\cmidrule{2-13}          & \multicolumn{1}{r}{EAM} & \multicolumn{1}{r}{CGMP} & \multicolumn{1}{r}{EAM} & \multicolumn{1}{r}{CGMP} & \multicolumn{1}{r}{EAM} & \multicolumn{1}{r}{CGMP} & \multicolumn{1}{r}{EAM} & \multicolumn{1}{r}{CGMP} & \multicolumn{1}{r}{EAM} & \multicolumn{1}{r}{CGMP} & \multicolumn{1}{r}{\textsc{EAM}} & \multicolumn{1}{r}{CGMP} \\
			\midrule
			High  & 1000  & 1000  & 1000  & 1000  & 1000  & 1000  & 468   & 304   & 922   & 1000  & 1000  & 1000 \\
			Mid   & 1000  & 1000  & 993   & 1000  & 1000  & 1000  & 108   & 32    & 688   & 720   & 1000  & 1000 \\
			Low   & 1000  & 1000  & 938   & 1000  & 1000  & 1000  & 0     & 0     & 405   & 332   & 1000  & 1000 \\
			\bottomrule
		\end{tabular}%
		\label{tab_comptime}					
	}
\end{table}

\section{Managerial insights}\label{insight}
The smooth operations of airlines require a good match between different types of resources, with the scarcest resource serving as the bottleneck of the operation. Consequently, one of the most frequently asked questions is, "What is the scarcity level of each resource?" The response to this question is important, as it will guide the strategic decisions for the future development of airlines.

However, to the best of our knowledge, there is no quantitative method in both academia and industry. For example, in our interview with a well-known worldwide airline, both the B737 and A320 crew members thought of themselves as the scarcest group, requiring immediate replenishment. However, the managers of the airline are also unclear about which group is more critical to operations. A similar situation holds for different fleet types. Therefore, in this study, we present a systematic way to analyze the scarcities of different resources based on the information from the solution to the TFACPP.

We first introduce the crew scarcity evaluating procedure and its application, followed by the counterparts for aircraft resources. Then, we analyze the matching degree between the crew and aircraft resources and provide tactical airline development strategies.

\subsection{Crew scarcity evaluating procedure and application}\label{crewvalue}
In this section, we provide a quantitative method, based on shadow prices of crew flight time constraints, to evaluate the scarcities of different types of crew resources.

In the CGMP, the $ b $th constraint in Eq. (\ref{mp_year}) represents a limitation for using at most $t_b$ hours of the crew flight time of fleet family $ b $. Thus, the corresponding dual variable ${\beta _b}$ can be regarded as the marginal profits of having one more hour of the crew flight time of $ b $ for the year. Then, intuitively, the yearly crew marginal profits are $\beta _b = {\beta _b} \times \bar t_b$. If $\beta _b$ is greater than zero, the crew flight time resources of $ b $ are exhausted; otherwise, the resources are surplus. The scarcer the crew resources of $ b $ are, the greater $\beta _b$ is.

\begin{figure}[htbp]
	\centering	
	\includegraphics[scale=.7]{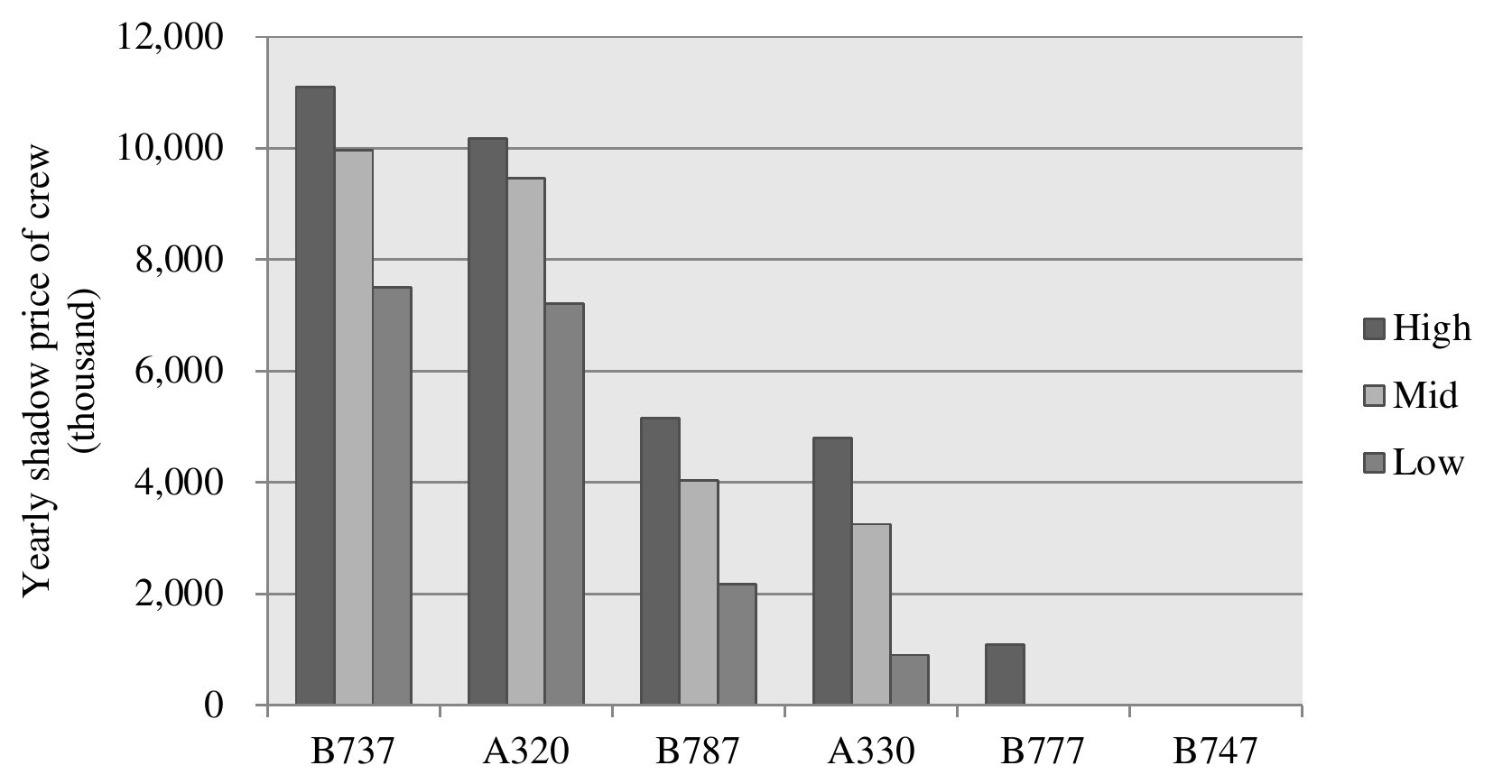}
	\caption{\centering{The yearly crew marginal profits of different fleet families and demand levels}}
	\label{fig_sdpcrew}
\end{figure}

Therefore, the yearly crew marginal profits provide key information regarding the scarcities of different types of crew resources, which can assist in crew hiring and training decisions. For example, the crew members in the fleet family with the largest $\beta _{{b^*}}$ are the scarcest and suggested to be prioritized in hiring and training. Moreover, airlines may consider transferring certain crew members from the less scarce type to the scarcer type through transition training. For the test cases described in Section \ref{caseintro}, we present the yearly crew marginal profits for all fleet families and demand levels in Fig. \ref{fig_sdpcrew}.

The trend in Fig. \ref{fig_sdpcrew} shows that, at all the demand levels, B737 and A320 correspond to the scarcest fleet families, followed by B787 and A320, and the resources for B747 and B777 may be redundant. Particularly, in the mid-demand level, it is recommended that the airline hires and trains crew members of B737, A320, B787, and A330, sequentially. Moreover, as there are surpluses in the crew flight time of B747 and B777, the airline may consider transferring crew members from B747 and B777 to B737 or postpone the crew transition from B737 to B747 and B777.

\subsection{Aircraft scarcity evaluating procedure and application}\label{aftvalue}
In this section, we evaluate the scarcities of different types of aircraft resources through the dual variables of the aircraft number constraints, so as to provide references for the above-mentioned decisions.

In each CGSP, the $ f $th constraint in Eq. (\ref{sp_afn}) represents a limitation for using at most ${k_f}$ number of aircraft of fleet type $ f $. Thus, the dual variable associated with this constraint, defined as $\gamma _f^m$, can be regarded as the marginal profits of having one more aircraft of fleet type $ f $ for month $ m $. Then, intuitively, the yearly aircraft marginal profits are $\gamma _f = \mathop \sum \limits_{m \in M} \gamma _f^m$. The fleet type with a positive $\gamma _f$ is likely worth expanding, as a larger $\gamma _f$ corresponds to a higher scarcity of the aircraft of fleet type $ f $.

\begin{figure}[h]
	\centering	
	\includegraphics[height=10cm,width=16cm]{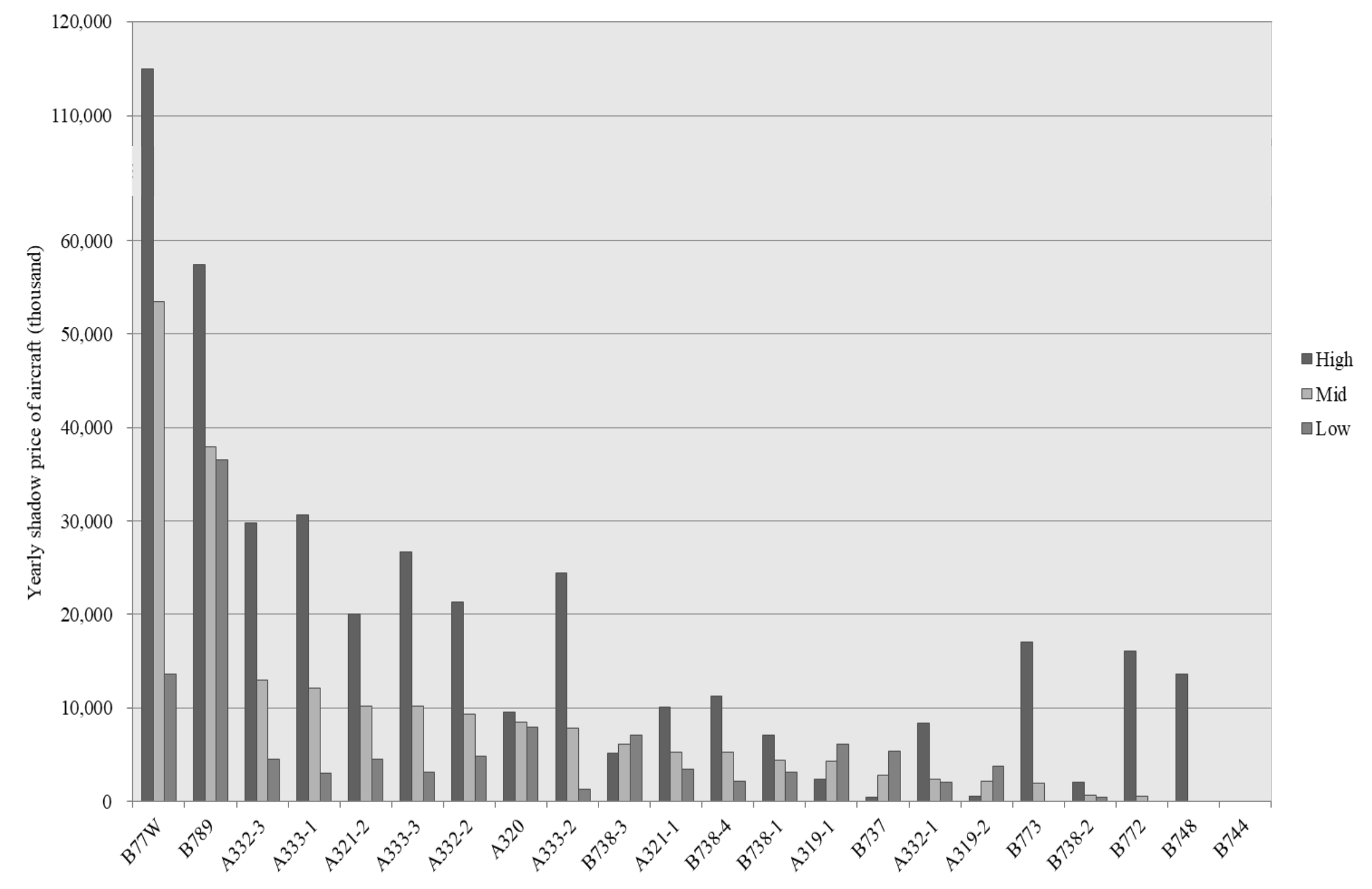}
	\caption{\centering{The yearly aircraft marginal profits of different fleet families and demand levels}}		\label{fig_sdpaft}
\end{figure}

Therefore, the yearly aircraft marginal profits reflect the scarcities of different types of aircraft resources, which can facilitate fleet planning decisions. For example, the aircraft of the fleet type with the largest $\gamma _{{f^*}}$ are the scarcest and suggested to be acquired first. In addition, $\gamma _f$ reflects the profitability of the aircraft of fleet type $ f $, which further affects the acquisition strategy (such as purchase and lease). In contrast, the fleet type with a small $\gamma _f$ may be redundant. Considering the test cases, we present the yearly aircraft marginal profits for all fleet types and demand levels in Fig. \ref{fig_sdpaft}.

As shown in Fig. \ref{fig_sdpaft}, certain similarities exist in the yearly aircraft marginal profits under different demand levels. For instance, fleet type B77W has the highest scarcity and the second one is fleet type B787 across all demand levels. Fleet type B744 has the lowest scarcity. Nevertheless, several differences can be observed among the demand levels. For example, in the high- and mid-demand levels, B77W is scarcer than B789; however, the opposite trend is noted for the low-demand level. Moreover, the scarcities of fleet types, such as A333-2, A319-1, B773, and B748, vary dramatically with the demand level. This finding highlights that the requirements for the aircraft resources of different fleet types vary with the demand levels.

\subsection{Matching degree evaluating procedure and application}\label{matchdegree}
Apart from the scarcities of crew and aircraft resources, the matching degree between crew and aircraft can also provide key managerial insights. In particular, owing to the variety of aircraft resources, it is necessary to obtain additional information to facilitate aircraft acquisition and replacement decisions. Moreover, as mentioned previously, each fleet type belongs to a fleet family for crew. Thus, the matching degree between $ f $ and $ b $ ($f \in {F_b}$) will provide rich insights for airline management. 

After obtaining the yearly marginal profits, we measure the matching degree between each fleet type and corresponding crew type (fleet family) by considering the similarity of the scarcities. Specifically, if both the fleet type and crew type have large yearly marginal profits and are reasonably scarce, they are well matched and worth expanding. If only one of them has large yearly marginal profits, they are not well matched and only the scarce one is worth expanding. They may be redundant if none of them has large yearly marginal profits.

\begin{figure}[h]
	\centering			
	\includegraphics[scale=0.5]{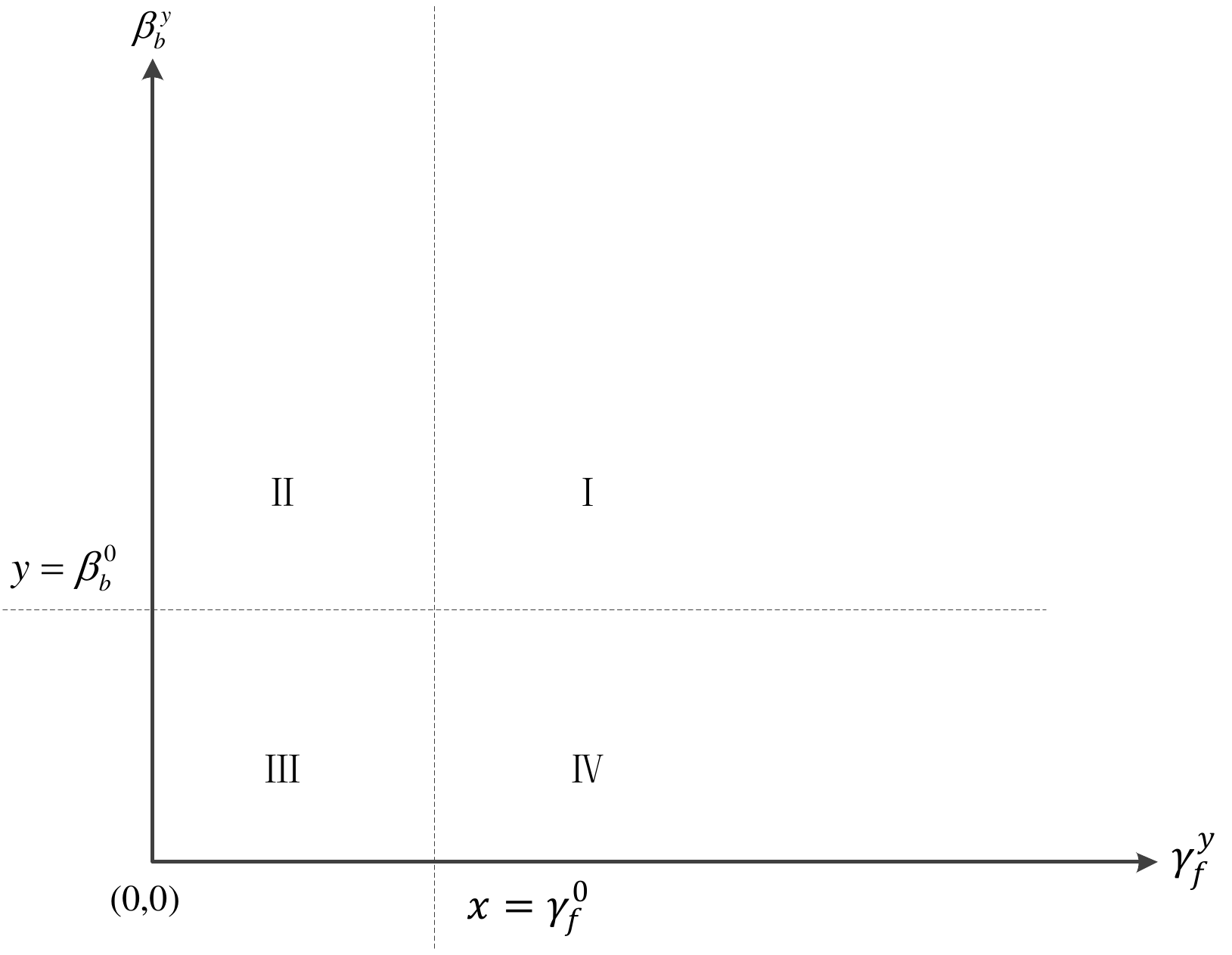}	
	\caption{\centering{General grouping of fleet and crew types}}
	\label{fig_generalGrp}
\end{figure}

As shown in Fig. \ref{fig_generalGrp}, we establish a two-dimensional coordinate system with $\gamma _f$ and $\beta _b$ as the $ x $- and $ y $-coordinates, respectively, and divide this coordinate system into four quadrants through $x = \gamma _f^0$ and $y = \beta _b^0$. We assume that the fleet type and crew type with yearly marginal profits larger than $\gamma _f^0$ and $\beta _b^0$, respectively, are reasonably scarce. In this case, for each $ f $ and $ b $ ($f \in {F_b}$), the quadrant that point ($\gamma _f,\;\beta _b$) lies in reflects the matching degree between them. The main characteristics of the quadrants are presented in Table \ref{tab_gengrp}.

\begin{table}[htbp]
	\caption{Characteristics of four quadrants}
	\centering
	{\def\arraystretch{1}  
		\begin{tabular}{lllp{8.5cm}}
			\toprule
			Quadrant & Crew scarcity & Aircraft scarcity & Matching degree between fleet and crew types \\
			\midrule
			\uppercase\expandafter{\romannumeral1}     & High  & High  & The crew and aircraft resources are well matched, and both resources are worth expanding. Consequently, the airline might develop the crew and aircraft resources in this quadrant. \\
			\uppercase\expandafter{\romannumeral2}     & High  & Low   & There is a lacking of crew resources in this quadrant, which might be worth expanding. \\
			\uppercase\expandafter{\romannumeral3}     & Low   & Low   & Neither resource is worth expanding. Consequently, the airline might retire the fleet type in this quadrant. The crew might be transferred to quadrant \uppercase\expandafter{\romannumeral1} or \uppercase\expandafter{\romannumeral2} for better use. \\
			\uppercase\expandafter{\romannumeral4}     & Low   & High  & There is a lacking of aircraft resources in this quadrant, which might be worth expanding. \\
			\bottomrule
		\end{tabular}%
		\label{tab_gengrp}
	}
\end{table}

\begin{figure}[h]
	\centering			
	\includegraphics[width=16cm,height=6cm]{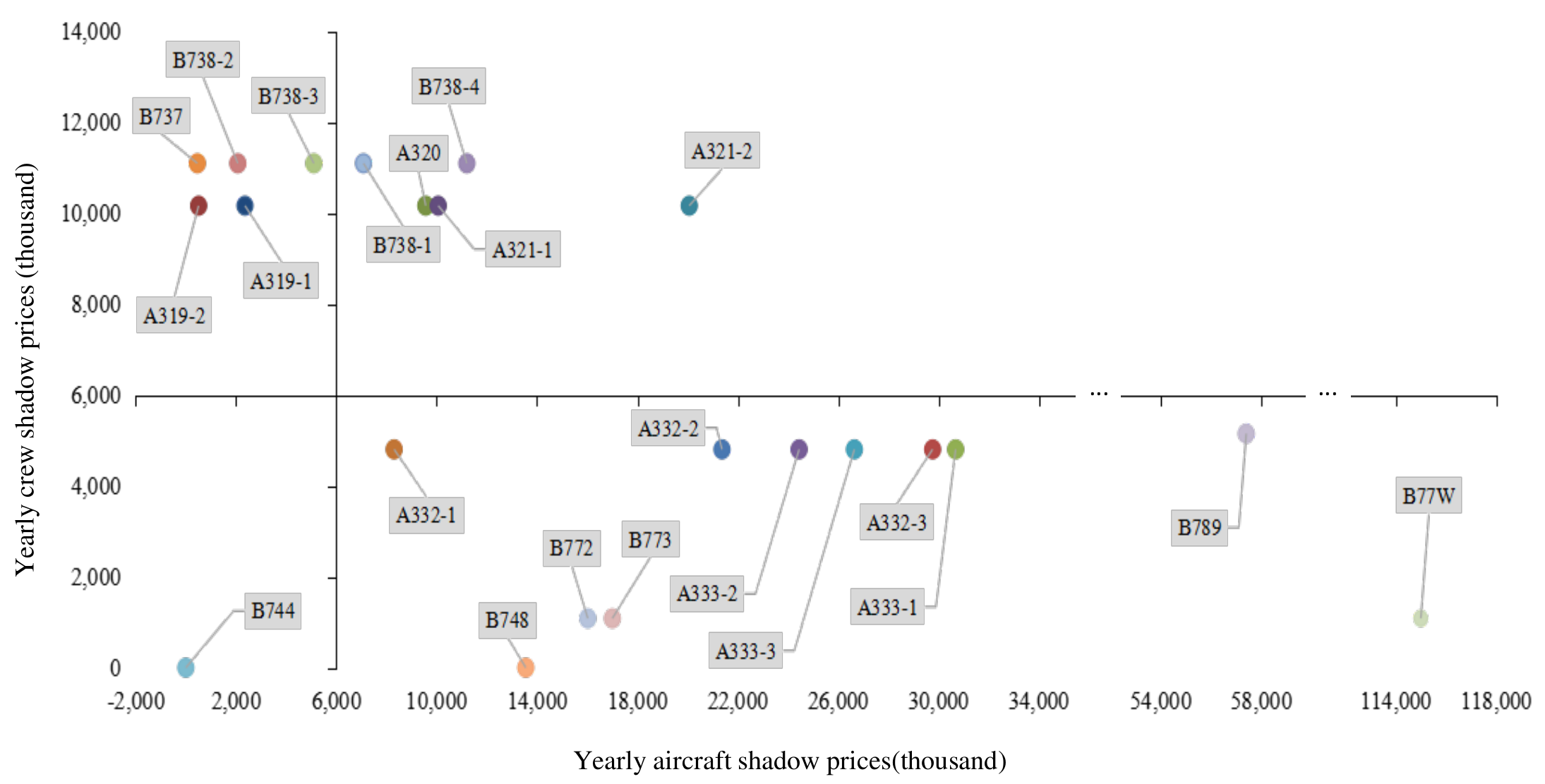}	
	\caption{\centering{Grouping of fleet types and crew types (High)}}
	\label{fig_highgrp}
\end{figure}

\begin{figure}[h]
	\centering			
	\includegraphics[width=16cm,height=6cm]{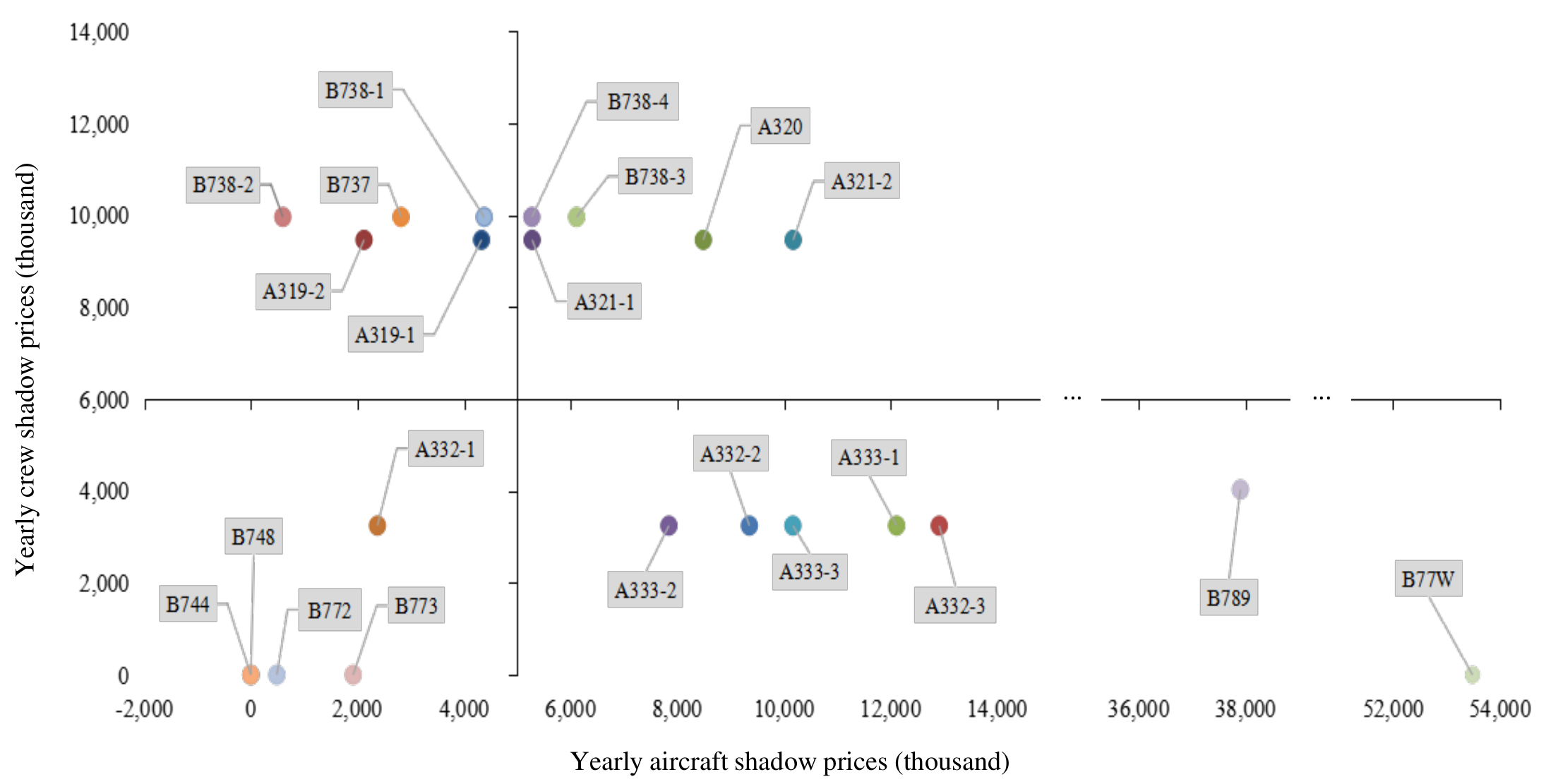}	
	\caption{\centering{Grouping of fleet types and crew types (Mid)}}
	\label{fig_midgrp}
\end{figure}

\begin{figure}[h]
	\centering			
	\includegraphics[width=16cm,height=6cm]{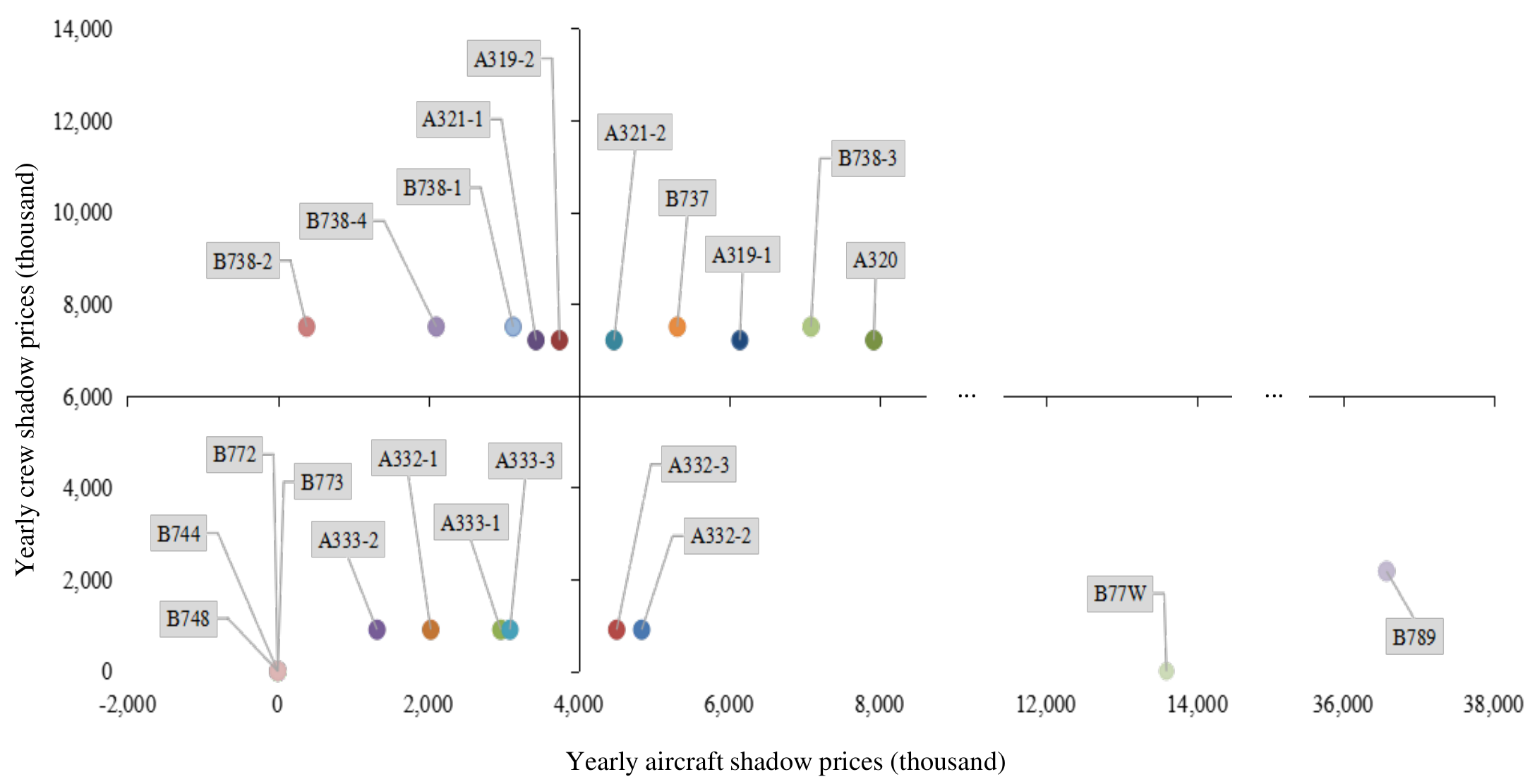}	
	\caption{\centering{Grouping of fleet types and crew types (Low)}}
	\label{fig_lowgrp}
\end{figure}

\begin{figure}[h]
	\centering			
	\includegraphics[width=13cm,height=6cm]{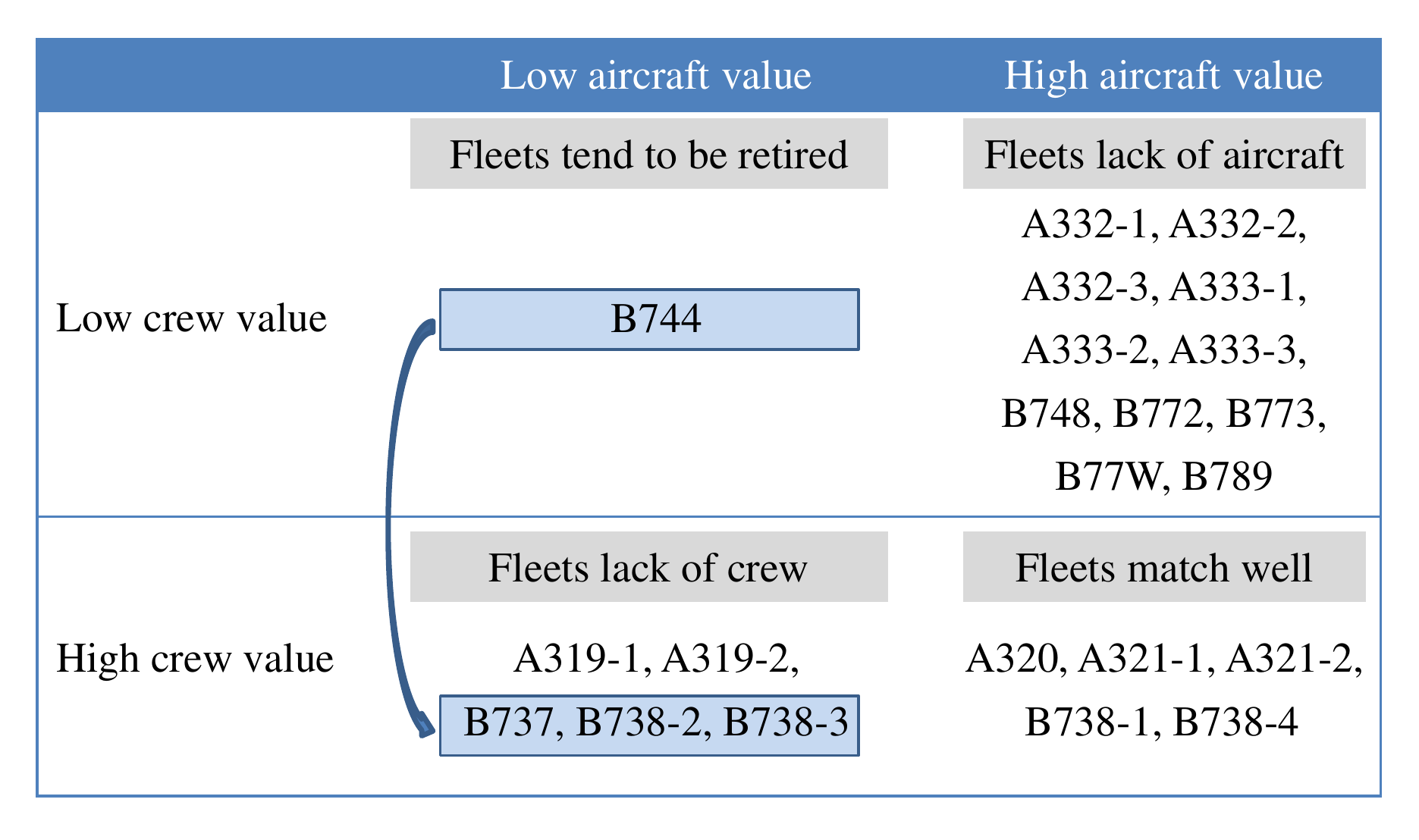}	
	\caption{\centering{Summary of grouping (High)}}
	\label{fig_highgrpsum}
\end{figure}

\begin{figure}[h]
	\centering			
	\includegraphics[width=13cm,height=6cm]{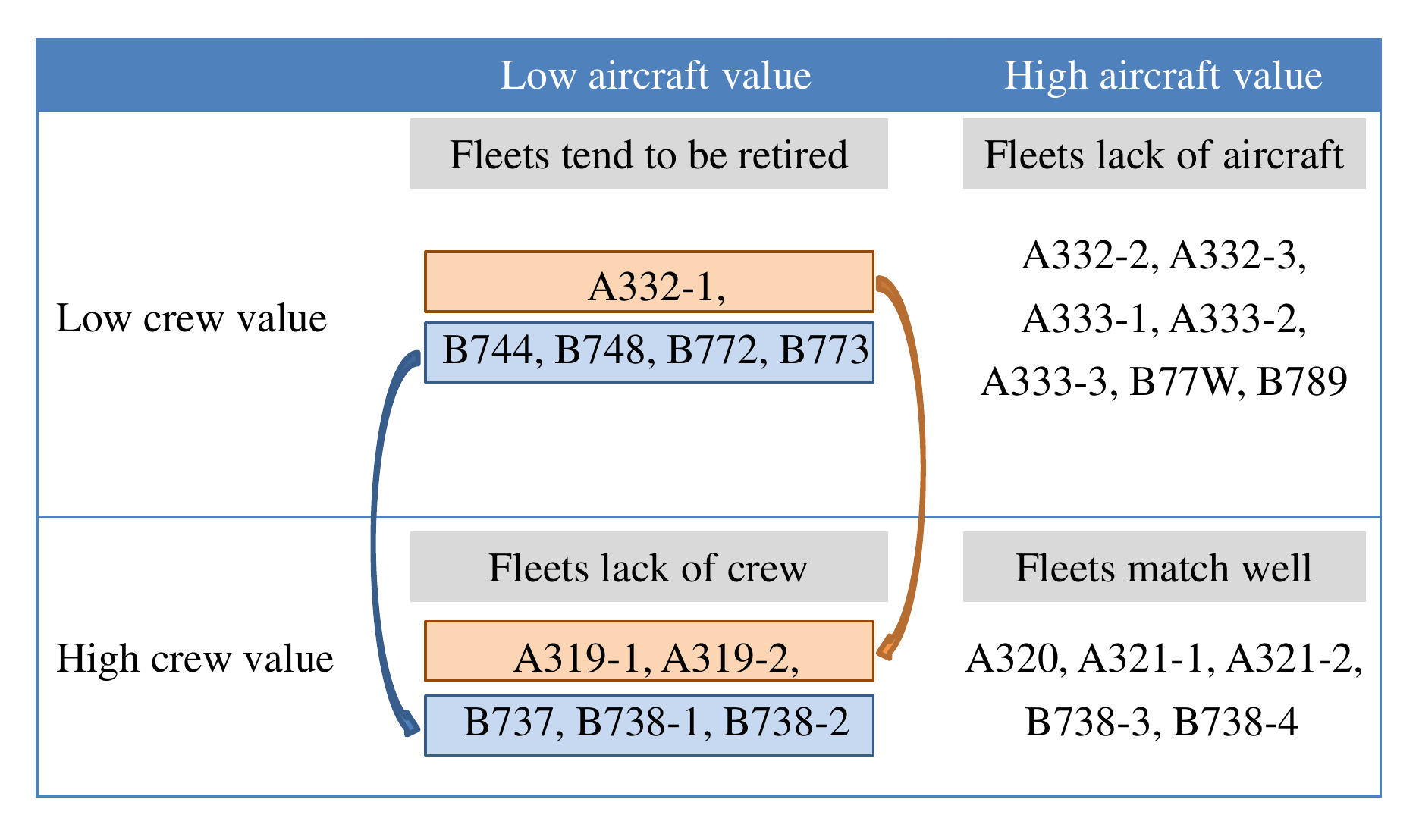}	
	\caption{\centering{Summary of grouping (Mid)}}
	\label{fig_midgrpsum}
\end{figure}

\begin{figure}[h]
	\centering			
	\includegraphics[width=13cm,height=6cm]{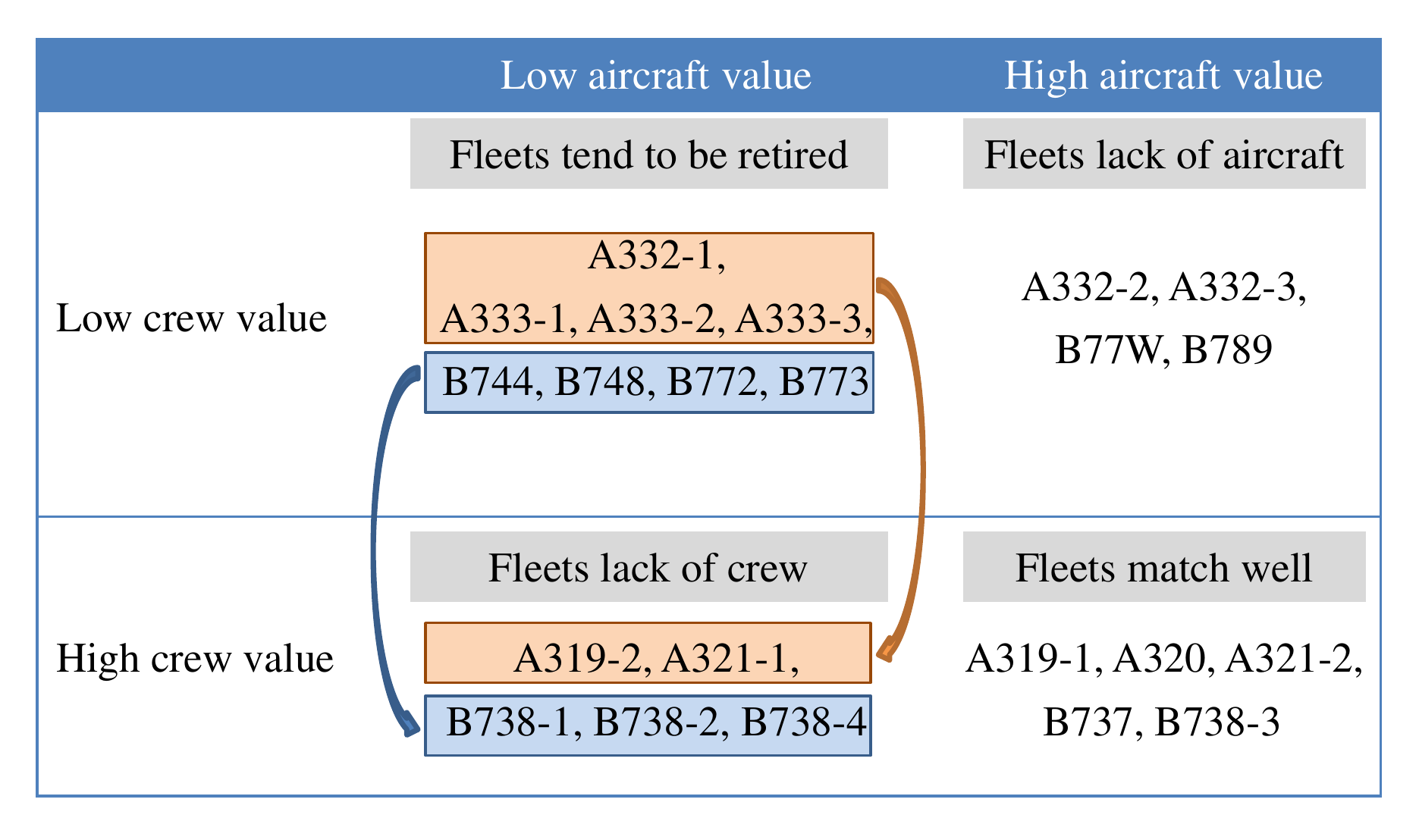}	
	\caption{\centering{Summary of grouping (Low)}}
	\label{fig_lowgrpsum}
\end{figure}

Finally, we consider the matching degree between each fleet type and corresponding crew type for the test cases with three demand levels. Based on the management experience, we set the thresholds of the yearly aircraft shadow prices $\gamma _f^0$ as 6, 5, and 4 million for the high-, mid-, and low-demand levels. Moreover, we set the threshold of the yearly crew shadow prices $\beta _b^0$ as 6 million for all three demand levels. The groupings of all fleet types and crew types according to the matching degree are shown in Fig. \ref{fig_highgrp}–\ref{fig_lowgrp} and summarized in Table \ref{tab_highgrp}–\ref{tab_lowgrp}. Particularly, we indicate a possible relocating method for crew members of the retired aircraft in Fig. \ref{fig_highgrpsum}–\ref{fig_lowgrpsum}.

\begin{table}[h]
	\caption{Grouping of fleet types and crew types (High)}
	\centering
	{\def\arraystretch{1}  
		\begin{tabular}{lllp{7cm}}
			\toprule
			Matching degree & \multicolumn{1}{l}{Aircraft scarcity} & \multicolumn{1}{l}{Crew scarcity} & Fleet types \\
			\midrule
			Well  & High  & High  & A320, A321-1, A321-2, B738-1, B738-4 \\
			Poor  & Low   & High  & A319-1, A319-2, B737, B738-2, B738-3 \\
			-     & Low   & Low   & B744 \\
			Poor  & High  & Low   & A332-1, A332-2, A332-3, A333-1, A333-2, A333-3, B748, B772, B773, B77W, B789 \\
			\bottomrule
		\end{tabular}%
		\label{tab_highgrp}						
	}
\end{table}

\begin{table}[h]
	\caption{Grouping of fleet types and crew types (Mid)}
	\centering
	{\def\arraystretch{1}  
		\begin{tabular}{lllp{7cm}}
			\toprule
			Matching degree & \multicolumn{1}{l}{Aircraft scarcity} & \multicolumn{1}{l}{Crew scarcity} & Fleet types \\
			\midrule
			Well  & High  & High  & A320, A321-1, A321-2, B738-3, B738-4 \\
			Poor  & Low   & High  & A319-1, A319-2,  B737, B738-1, B738-2 \\
			-     & Low   & Low   & A332-1, B744, B748, B772, B773 \\
			Poor  & High  & Low   & A332-2, A332-3, A333-1, A333-2, A333-3, B77W, B789 \\
			\bottomrule
		\end{tabular}%
		\label{tab_midgrp}
	}
\end{table}

\begin{table}[h]
	\caption{Grouping of fleet types and crew types (Low)}
	\centering
	{\def\arraystretch{1}  
		\begin{tabular}{lllp{7cm}}
			\toprule
			Matching degree & \multicolumn{1}{l}{Aircraft scarcity} & \multicolumn{1}{l}{Crew scarcity} & Fleet types \\
			\midrule
			Well  & High  & High  & A319-1, A320, A321-2, B737, B738-3 \\
			Poor  & Low   & High  & A319-2, A321-1, B738-1, B738-2, B738-4 \\
			-     & Low   & Low   & A332-1, A333-1, A333-2, A333-3, B744, B748, B772, B773 \\
			Poor  & High  & Low   & A332-2, A332-3, B77W, B789 \\
			\bottomrule
		\end{tabular}%
		\label{tab_lowgrp}
	}
\end{table}

The matching degree can provide several managerial insights, which are meaningful to airlines:

\begin{enumerate}
	\item The crew and aircraft resources corresponding to B744 may need to be reduced. In all the demand levels, the crew and aircraft resources of B744 appear to be surplus. Especially, all the associated shadow prices (fleet types B744 and B748) are zero at the mid- and low-demand levels.
	
	\item The crew and aircraft resources corresponding to fleet types A320 and A321-2 may be worth expanding. In all the demand levels, the crew and aircraft resources related to A320 and A321-2 match well and are both scarce. In addition, fleet types A321-2, B738-3, and B738-4 may also be worth expanding, at least in two of three demand levels.
	
	\item The aircraft resources corresponding to fleet type B77W may be expanded. In all the demand levels, the related crew resources (B777) are redundant; however, the aircraft resources are highly scarce. Moreover, the fleet type B789 may be worth expanding, although the scarcity of the related crew resources (crew type B787) is mid-level.
	
	\item The crew resources of crew types B737 and A320 may be worth expanding. In all the demand levels, several fleet types pertaining to B737 and A320 are scarce in crew resources but redundant in aircraft resources. These findings can also be derived through the yearly crew shadow prices.
\end{enumerate}

It must be noted that in the two years spent conducting this research, the airline retired the aircraft of the fleet type B744. Moreover, the total number of aircraft of the airline increased by 34, with the number of aircraft for fleet types B77W and B789 increasing by 5 and 7, respectively. These observations show that the results of this study are broadly consistent with managerial experience. Moreover, our study can provide a quantitative reference for managers.

\section{Conclusions}\label{conclu}
Motivated by the mismatch between the crew and aircraft resources in the aviation industry, we address the TFACPP and reformulate it by the Benders decomposition. The major difficulty of the problem lies in the BMP, and we develop a column generation algorithm to solve it efficiently. We evaluate our method using a realistic test case with three demand levels. The computational results demonstrate the effectiveness of the solution method and the high quality of its solutions. For all the demand levels, the computation time is acceptable and the integrality gaps between the MIP and LP solutions are very small. Moreover, compared with the solutions to the EAM, the TFATAP can increase the profit by more than two hundred million CNY. 

In addition, we provide a quantitative method to measure the scarcities of crew and aircraft resources and the matching degree between crew and aircraft based on the shadow prices corresponding to relative constraints. These information can provide useful managerial insights regarding the acquisition, replacement, and transition of the crew and aircraft resources. For example, the findings demonstrate the significant redundancy in fleet family B747 and fleet type B744 and the significant shortage of aircraft resources for fleet types B77W and B789. The actual trend of the airline indicates that the proposed managerial insights are consistent with the managerial experience while being more scientific and accurate.

\ACKNOWLEDGMENT{%
	This study is supported by National Natural Science Foundation of China under Grant No.71825001.
}

%
%
%


\bibliographystyle{informs2014trsc} 
\bibliography{mybibfile} 


\end{document}